\documentclass[a4paper,11pt]{article}
\usepackage{mfpaperstuff2}
\usepackage{mfabacus}
\usetikzlibrary{calc}

\begin{document}

\newcommand\lb[2]{#1^{(#2)}}
\newcommand\reg{^{\operatorname{reg}}}
\newcommand\ws[2]{W^{#1,#2}}
\newcommand\sws[2]{\spin W^{#1,#2}}
\newcommand\pw[2]{\calp^{#1,#2}}
\newcommand\spw[2]{\spin\calp^{#1,#2}}
\newcommand\sspe[1]{\mathbf{S}^{#1}}
\newcommand\sid[1]{\mathbf{D}^{#1}}
\newcommand\infi{\raisebox{1pt}{$\infty$}}
\newcommand\lr[3]{\operatorname{c}^{#1}_{#2#3}}
\newcommand\hf\iota
\newcommand\fh{\bar\iota}
\newcommand\cls\preccurlyeq
\newcommand\ca[1]{\stackrel{#1}\leadsto}
\newcommand\car{\ca r}
\newcommand\cak{\ca k}
\newcommand\znz{\bbz/n\bbz}
\newcommand\zmz{\bbz/m\bbz}
\newcommand\fonicel{\spin\fo^l}
\newcommand\f[1]{\mathrm f_{#1}}
\renewcommand\t[1]{\mathrm t_{#1}}
\newcommand\spf[1]{\spin{\mathrm{f}}_{#1}}
\newcommand\e[1]{\mathrm e_{#1}}
\newcommand\spe[1]{\spin{\mathrm{e}}_{#1}}
\newcommand\spt[1]{\spin{\mathrm{t}}_{#1}}
\newcommand\jo[2]{#2|#1}
\newcommand\cb[1]{\mathrm{G}_{#1}}
\newcommand\scb[1]{\spin{\mathrm{G}}_{#1}}
\renewcommand\labelenumi{(\theenumi)}
\newcommand\tac\textasteriskcentered
\newcommand\asa\Theta
\newcommand\sas{\asa^{-1}}
\newcommand\bsspu[1]{\operatorname{BSSP}_{#1}}
\newcommand\bssp{\bsspu{}}
\newcommand\sspu[1]{\operatorname{SSP}_{#1}}
\newcommand\ssp{\sspu{}}
\newcommand\skw[2]{#1\setminus#2}
\newcommand\hs[2]{\operatorname{hs}_{#1#2}(t)}
\newcommand\hsb[2]{\widebar{\operatorname{hs}}_{#1#2}}
\newcommand\ddd[1]{\partial_{#1}}
\newcommand\ar{\operatorname{ar}}
\newcommand\jh[2]{#2\|#1}
\newcommand\dn[2]{\operatorname{d}_{#1#2}}
\newcommand\sdn[2]{\spin{\operatorname{d}}_{#1#2}}
\newcommand\miabb{\raisebox{-1.25pt}{$\begin{tikzpicture}[scale=1.1]\fill[gray!30!white,rounded corners=1.5pt](-.15,-.18)rectangle++(.6,.36);\draw(0,-.16)--++(0,.32);\draw(0.3,-.16)--++(0,.32);\shade[shading=ball,ball color=black](0,0)circle(.13);\shade[shading=ball,ball color=black](0.3,0)circle(.13);\end{tikzpicture}$}}
\newcommand\miabn{\raisebox{-1.25pt}{$\begin{tikzpicture}[scale=1.1]\fill[gray!30!white,rounded corners=1.5pt](-.15,-.18)rectangle++(.6,.36);\draw(.24,0)--++(.12,0);\draw(0,-.16)--++(0,.32);\draw(0.3,-.16)--++(0,.32);\shade[shading=ball,ball color=black](0,0)circle(.13);\end{tikzpicture}$}}
\newcommand\mianb{\raisebox{-1.25pt}{$\begin{tikzpicture}[scale=1.1]\fill[gray!30!white,rounded corners=1.5pt](-.15,-.18)rectangle++(.6,.36);\draw(-.06,0)--++(.12,0);\draw(0,-.16)--++(0,.32);\draw(0.3,-.16)--++(0,.32);\shade[shading=ball,ball color=black](0.3,0)circle(.13);\end{tikzpicture}$}}
\newcommand\mibbb{\raisebox{-1.25pt}{$\begin{tikzpicture}[scale=1.1]\fill[gray!30!white,rounded corners=1.5pt](-.15,-.18)rectangle++(.6,.36);\draw(0,-.16)--++(0,.32);\draw(0.3,-.16)--++(0,.32);\shade[shading=ball,ball color=white](0,0)circle(.13);\shade[shading=ball,ball color=white](0.3,0)circle(.13);\end{tikzpicture}$}}
\newcommand\mibbn{\raisebox{-1.25pt}{$\begin{tikzpicture}[scale=1.1]\fill[gray!30!white,rounded corners=1.5pt](-.15,-.18)rectangle++(.6,.36);\draw(.24,0)--++(.12,0);\draw(0,-.16)--++(0,.32);\draw(0.3,-.16)--++(0,.32);\shade[shading=ball,ball color=white](0,0)circle(.13);\end{tikzpicture}$}}
\newcommand\mibnn{\raisebox{-1.25pt}{$\begin{tikzpicture}[scale=1.1]\fill[gray!30!white,rounded corners=1.5pt](-.15,-.18)rectangle++(.6,.36);\draw(.24,0)--++(.12,0);\draw(-.06,0)--++(.12,0);\draw(0,-.16)--++(0,.32);\draw(0.3,-.16)--++(0,.32);\end{tikzpicture}$}}
\newcommand\mibnb{\raisebox{-1.25pt}{$\begin{tikzpicture}[scale=1.1]\fill[gray!30!white,rounded corners=1.5pt](-.15,-.18)rectangle++(.6,.36);\draw(-.06,0)--++(.12,0);\draw(0,-.16)--++(0,.32);\draw(0.3,-.16)--++(0,.32);\shade[shading=ball,ball color=white](0.3,0)circle(.13);\end{tikzpicture}$}}
\newcommand\sep{separated\xspace}
\newcommand\ssep{super-\sep}
\newcommand\bsep{bar-separated\xspace}
\newcommand\sbsep{super-\bsep}
\newcommand\ber[2]{\calb_{#1}(#2)}
\newcommand\abdisp{abacus display}
\newcommand\babdisp{bar-abacus display}
\newcommand\babd{\babdisp\xspace}
\newcommand\babds{\babdisp s\xspace}
\newcommand\ab[1]{\operatorname{Ab}(#1)}
\newcommand\bab[1]{\widebar{\operatorname{Ab}}(#1)}
\newcommand\pl[1]{\calp_{\ls #1}}
\newcommand\pml{\calp^+_{\ls l}}
\newcommand\pnice[1]{\calp_{#1}^{\operatorname{std}}}
\newcommand\ip[2]{\left(#1,#2\right)}
\newcommand\ipx[2]{(#1,#2)}
\newcommand\hi{\hat\imath}
\newcommand\aarr[1]{\stackrel{#1}\Longrightarrow}
\newcommand\arr[1]{\stackrel{#1}\longrightarrow}
\newcommand\len[1]{l(#1)}
\newcommand\ol[1]{\widebar{#1}}
\newcommand\ba{$h$-strict}
\newcommand\hstr{\calp^{(h)}}
\newcommand\fo{\calf}
\newcommand\nice{standard\xspace}
\newcommand\bpx{$h$-strict partition}
\newcommand\bp{\bpx\xspace}
\newcommand\bps{\bpx s\xspace}
\newcommand\spin{\check}
\newcommand\hsss{\hat{\mathfrak{S}}_}

\title{Comparing Fock spaces in types $A^{(1)}$ and $A^{(2)}$}


\toptitle

\begin{abstract}
We compare the canonical bases of level-$1$ quantised Fock spaces in affine types $A^{(1)}$ and $A^{(2)}$, showing how to derive the canonical basis in type $A^{(2)}_{2n}$ from the the canonical basis in type $A^{(1)}_n$ in certain weight spaces. In particular, we derive an explicit formula for the canonical basis in extremal weight spaces, which correspond to RoCK blocks of double covers of symmetric groups. In a forthcoming paper with Kleshchev and Morotti we will use this formula to find the decomposition numbers for RoCK blocks of double covers with abelian defect.
\end{abstract}

\tableofcontents

\section{Introduction}

This paper is motivated by the \emph{decomposition number problem} for the symmetric groups and their double covers in characteristic $p$. Although a solution to this problem seems a long way off, several important results are known. One of these results gives the decomposition numbers for \emph{RoCK blocks} of symmetric groups; these are particularly well understood blocks which have been used in a variety of applications. The formula for the decomposition numbers for RoCK blocks in the abelian defect case was given by Chuang and Tan \cite{cts}, and results of Turner \cite{turn} allow these results to be extended to RoCK blocks with non-abelian defect groups.

It is natural to seek analogous results for the double cover $\hsss n$ of the symmetric group (which controls projective representations of $\sss n$). The case of characteristic $2$ behaves very differently (and is dealt with in \cite{mfspin2alt}), so we concentrate here on odd characteristic. The representations of $\hsss n$ which do not descend to representations of $\sss n$ are called \emph{spin} representations of $\sss n$, and the blocks of $\hsss n$ containing spin representations are called spin blocks.  RoCK blocks for symmetric groups can be characterised as elements of the maximal equivalence class of blocks under the Scopes equivalence \cite{sco} on blocks of symmetric groups. The Scopes--Kessar equivalence \cite{kes} for spin blocks of double covers suggests a natural analogue of RoCK blocks, and this can be realised in a combinatorial way using the abacus. These blocks have recently been studied in detail by Kleshchev and Livesey \cite{kl}, who prove Brou\'e's abelian defect group conjecture for RoCK blocks. This has been used even more recently by Brundan and Kleshchev, and independently Ebert, Lauda and Vera, to show that Brou\'e's conjecture holds for all spin blocks \cite{bk,elv}. However, the results of Kleshchev and Livesey do not directly address the decomposition number problem for spin RoCK blocks, and this is the main focus here.

In this paper we address RoCK blocks by studying quantum algebra. Let $U=U_q(A^{(1)}_{p-1})$ denote the quantised universal enveloping algebra of the affine Kac--Moody algebra of type $A^{(1)}_{p-1}$. The \emph{level-$1$ Fock space} is a highest-weight $U$-module with a simple combinatorial construction in terms of integer partitions. The submodule generated by a highest-weight vector is isomorphic to the irreducible highest-weight module $V(\La_0)$, so the Fock space provides a combinatorial framework for studying $V(\La_0)$; this approach has proved useful, for example, in constructing the crystals. $V(\La_0)$ possess an important basis called the \emph{canonical basis}, which provides a connection to representation theory of symmetric groups and Iwahori--Hecke algebras, via the work of Lascoux--Leclerc--Thibon \cite{llt} and Ariki \cite{ari}, who showed that decomposition numbers for Hecke algebras of type $A$ in characteristic zero can be obtained by specialising canonical basis coefficients at $q=1$. This means in particular that these decomposition numbers can be calculated algorithmically. A further conjecture due to James suggested that the same should apply for decomposition numbers of symmetric groups, in blocks with abelian defect groups. James's conjecture is now known to be false in general \cite{willi}, but there are a wide variety of situations where it is known to hold, and it has provided important inspiration for results on decomposition numbers. In particular, the formula due to Chuang and Tan \cite{cts} for decomposition numbers for RoCK blocks of $\sss n$ was inspired by their earlier calculation of the canonical basis in weight spaces corresponding to RoCK blocks, and shows in particular that James's conjecture holds for RoCK blocks.

An analogous connection to quantum groups for spin representations was found by Leclerc and Thibon \cite{lt}, using the quantum group of type $A^{(2)}_{p-1}$. This quantum group also acts on a combinatorially defined level-$1$ Fock space (now defined in terms of $p$-strict partitions), which possesses an irreducible highest-weight submodule with a canonical basis. Leclerc and Thibon formulated an analogue of James's conjecture for spin representations of $\hsss n$. This conjecture is also known not to hold in general, but it does hold in many special cases, in particular all known cases of decomposition numbers for RoCK blocks. Motivated by this conjecture, the aim of the present paper is to find the canonical bases for the weight spaces of the basic $U_q(A^{(2)}_{p-1})$-module corresponding to spin RoCK blocks. Rather than directly determining the canonical basis, we deduce our result from the results of Chuang--Tan by proving more general results comparing the canonical bases in types $A^{(1)}$ and $A^{(2)}$: we show that if $\be$ is a restricted $p$-strict partition satisfying a particular additional condition which says that the $p$-bar-core of $\be$ is large in a certain specific sense relative to the sum of the parts of $\be$ divisible by $p$, then we can obtain the canonical basis element labelled by $\be$ from a corresponding canonical basis element in type $A^{(1)}$ by an adjustment involving inverse Kostka polynomials. In certain cases this allows the canonical bases for entire weight spaces to be computed, including weight spaces corresponding to RoCK blocks. Combining this result with the Chuang--Tan formula for RoCK blocks in type $A^{(1)}$ yields our main result (\cref{mainrouq}).

Combining our theorem with the Leclerc--Thibon conjecture (specialised to the case of RoCK blocks), we arrive at a conjecture for the decomposition numbers for RoCK blocks (\cref{rockconj}). In a forthcoming paper with Kleshchev and Morotti \cite{fkm} we will prove this conjecture.

\begin{ack}
This research was partly supported by EPSRC Small Grant EP/W005751/1.
\end{ack}

\section{Background}

In this section we set out some background details on partitions and Fock spaces.

\subsection{Elementary notation}

We write $\bbn$ for the set of positive integers and $\bbn_0=\bbn\cup\{0\}$. Given $m\in\bbn$, the set $\zmz$ consists of cosets $a+m\bbz=\lset{a+mb}{b\in\bbz}$. Given any set $B\subseteq\bbz$ and $a\in\bbz$, we write $B+a=\lset{b+a}{b\in B}$.

\subsection{Partitions}

A \emph{partition} is an infinite weakly decreasing sequence $\la=(\la_1,\la_2,\dots)$ of non-negative integers which is eventually zero. When writing partitions, we omit the trailing zeroes and group together equal parts with a superscript. The partition $(0,0,\dots)$ is written as $\vn$. If $\la$ is a partition, the integers $\la_1,\la_2,\dots$ are called the \emph{parts} of $\la$. We write $\card\la=\la_1+\la_2+\cdots$, and we say that $\la$ is a partition of $\card\la$. The \emph{length} $\len\la$ is the number of non-zero parts of $\la$. We write $\calp$ for the set of all partitions.

If $\la$ is a partition and $n\in\bbn$, then $n\la$ is defined to be the partition $(n\la_1,n\la_2,\dots)$. If $\la$ and $\mu$ are partitions, then $\la\sqcup\mu$ is defined to be the partition of $\card\la+\card\mu$ obtained by combining all the parts of $\la$ and $\mu$ in decreasing order.

The \emph{Young diagram} of $\la$ is the set
\[
\lset{(r,c)\in\bbn^2}{c\ls\la_r}
\]
whose elements are called the \emph{nodes} of $\la$. In general, a node means an element of $\bbn^2$. We use the English convention for drawing Young diagrams, in which $r$ increases down the page and $c$ increases from left to right. We abuse notation by identifying $\la$ with its Young diagram; so for example we may write $\la\subseteq\mu$ to mean that $\la_r\ls\mu_r$ for all $r$.

If $\la$ is a partition, the \emph{conjugate} partition $\la'$ is the partition obtained by reflecting the Young diagram of $\la$ in the main diagonal; that is, $\la'_r=\card{\lset{c\in\bbn}{\la_c\gs r}}$.

The \emph{dominance order} on partitions is defined by writing $\la\dom\mu$ (and saying that $\la$ \emph{dominates} $\mu$) if $\card\la=\card\mu$ and  $\la_1+\dots+\la_r\gs\mu_1+\dots+\mu_r$ for all $r$.

We say that a node $(r,c)$ of $\la$ is \emph{removable} if it can be removed from~$\la$ to leave the Young diagram of a partition (that is, if $c=\la_r>\la_{r+1}$), and we write the resulting partition as $\la\setminus(r,c)$. Similarly, a node $(r,c)$ not in $\la$ is an \emph{addable node} of $\la$ if it can be added to $\la$ to give a partition, and we write this partition as $\la\cup(r,c)$.

Now fix an integer $m\gs2$. We say that a partition $\la$ is \emph{$m$-restricted} if $\la_r-\la_{r+1}<m$ for all $r$. We define the \emph{residue} of the node $(r,c)$ to be $c-r+m\bbz$. Given $i\in\zmz$, we use the term \emph{$i$-addable node} to mean ``addable node of residue $i$'', and we define the term \emph{$i$-removable} similarly. If $\la$ and $\mu$ are partitions, we write $\la\arr{i:r}\mu$ to mean that $\mu$ is obtained from~$\la$ by adding $r$ nodes of residue $i$. (In the case $r=1$, we just write $\la\arr i\mu$.)

We also need to define rim hooks and $m$-cores. The \emph{rim} of a partition $\la$ is the set of nodes $(r,c)$ of $\la$ for which $(r+1,c+1)$ is not a node of $\la$. A \emph{rim $m$-hook} of $\la$ is a set of $m$ consecutive nodes of the rim which can be removed to leave the Young diagram of a smaller partition. The \emph{$m$-core} of $\la$ is the partition obtained by repeatedly removing rim $m$-hooks from $\la$ until none remain, and the \emph{$m$-weight} of $\la$ is the number of rim hooks removed to reach the $m$-core.

It is convenient in combinatorial representation theory to depict partitions using the abacus. Keeping $m$ fixed as above, we draw an abacus with $m$ vertical runners labelled $0,\dots,m-1$ from left to right. We mark positions on the runners, labelled $0,1,2,\dots$ from left to right along successive rows from top to bottom.

Now given a partition $\la$, we choose a large integer $s$, and define the \emph{beta-set}
\[
\ber s\la=\lset{\la_r+s-r}{1\ls r\ls s}.
\]
Now draw the \emph{$s$-bead \abd for $\la$} by placing beads on the abacus in the positions corresponding to all the elements of $\ber s\la$. In an \abd, we will say that a position is \emph{occupied} if there is a bead at that position, and \emph{empty} or \emph{unoccupied} otherwise.

For example, suppose $m=4$ and $\la=(8,7,5^2,2,1^3)$. Choosing $s=9$, we obtain the following \abd. (Whenever we draw an \abd, we adopt the convention that all positions below those depicted are empty.)
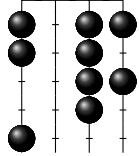
\[
\abacus(lmmr,bnbb,bnbn,nnbb,nnbn,bnnn)
\]

The \abd for a partition is very useful in two ways.
\begin{enumerate}
\item
Given $i\in\zmz$, let $a\in\{0,\dots,e-1\}$ be such that $i=a-s+m\bbz$. Then $i$-addable nodes of $\la$ correspond to empty positions $b$ on runner $a$ such that position $b-1$ is occupied (or $b=0$), while $i$-removable nodes correspond to occupied positions $b\gs1$ on runner $a$ for which position $b-1$ is empty.
\item
If position $b\gs m$ is occupied while position $b-m$ is empty, then moving the bead from position $b$ to position $b-m$ corresponds to removing a rim $m$-hook from $\la$. This means in particular that the \abd for the $m$-core of $\la$ can be obtained by moving all the beads up their runners as far as they will go, and the number of bead moves needed to do this is the $m$-weight of $\la$.
\end{enumerate}

We end with two combinatorial lemmas we will need later on. Given $\la,\mu\in\calp$ and $r\in\bbn_0$, we write $\la\car\mu$ to mean that $\mu$ is obtained from $\la$ by adding $r$ nodes in different columns. We want to interpret this condition in terms of beta-sets; the next \lcnamecref{carbeta} follows easily from the definitions.

\begin{lemma}\label{carbeta}
Suppose $\la,\mu\in\calp$ and $r\in\bbn_0$.
\begin{enumerate}
\item
$\la\car\mu$ \iff there is a set $A\subset\bbn$ with $\card A=r$ such that $\mu$ is obtained from $\la$ by adding a part equal to $a$ for each $a\in A$ and then removing a part equal to $a-1$ for each $a\in A$.
\item
Suppose $s\in\bbn$ is large. Then $\la\car\mu$ \iff there is a set $A\subset\bbn$ with $\card A=r$ and $A\cap\ber s\la=\emptyset$ such that $\ber s\mu=\ber s\la\cup A\setminus(A-1)$.
\end{enumerate}
\end{lemma}

We will also need the following lemma.

\begin{lemma}\label{carlem}
Suppose $\pi\in\calp$ with $\pi\neq\vn$. Let $k$ equal the last non-zero part of $\pi$, and let $\pi^-$ be the partition obtained by removing this last part. If $\si,\rho\in\calp$ with $\si\dom\pi^-$ and $\si\cak\rho$, then $\rho\dom\pi$, with equality only if $\si=\pi^-$.
\end{lemma}

\begin{pf}
Let $\si^+$ be the partition obtained from $\si$ by adding a node in each of the first $k$ columns. Since $\rho$ is obtained from $\si$ by adding nodes in some $k$ different columns, it follows that $\rho\dom\si^+$. Also, from the definition of the dominance order (and the fact that conjugation of partitions reverses the dominance order) the condition $\si\dom\pi^-$ implies that $\si^+\dom\pi$. So $\rho\dom\pi$, and in order to get equality we need $\si^+=\pi$, which is the same as $\si=\pi^-$.
\end{pf}

\subsection{\bps}

Now suppose $h\gs3$ is odd. We say that a partition is \emph{\ba} if there is no $r$ for which $\la_r=\la_{r+1}\nequiv0\ppmod h$. We write $\hstr$ for the set of all \bps. 
We will generally use letters near the start of the Greek alphabet for \bps. We define the \emph{bar-residue} of a node $(r,c)$ to be the smaller of the residues of $c-1$ and $h-c$ modulo $h$; so the bar-residue of a node depends only on the column in which it lies, and the bar-residues follow the pattern
\Yboxdim{13pt}
\[
\gyoungs(1.4,:0:1:2^2\cdots:{{\frac{h-1}2}}:{{\frac{h+1}2}}:{{\frac{h-1}2}}^2\cdots:2:1:0:0:1:2^2\cdots)
\]
from left to right. A node of an \bp $\la$ is \emph{bar-removable} if it can be removed from~$\la$, possibly together with some other nodes of the same bar-residue, to leave an \bp. Given $i\in\{0,\dots,\frac{h+1}2\}$, an \emph{$i$-bar-removable} node means a bar-removable node of bar-residue $i$. We define bar-addable and $i$-bar-addable nodes in a similar way. Note the distinction between removable nodes and bar-removable nodes. For example, if $h=5$ and $\alpha=(6,2,1)$, then $(2,2)$ is removable but not bar-removable, while $(1,5)$ is bar-removable but not removable. This can be seen from the following diagram, in which we label the nodes with their bar-residues.
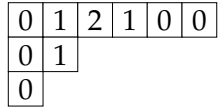
\[
\young(012100,01,0)
\]

If $\al$ and $\be$ are \bps, we write $\al\aarr{i:r}\be$ to mean that $\be$ is obtained from $\al$ by adding $r$ nodes of bar-residue $i$. (In the case $r=1$, we just write $\al\aarr i\be$.)

Now we define the analogue of $m$-cores. If $\al$ is an \bp, then removing an $h$-bar from $\la$ means either
\begin{itemize}
\item
replacing a part $\al_r\gs h$ with $\al_r-h$ (provided either $\al_r\equiv0\ppmod h$ or $\al_r-h$ is not a part of $\al$) and reordering, or
\item
deleting two parts which sum to $h$.
\end{itemize}
The \emph{$h$-bar-core} of $\al$ is the \bp obtained by repeatedly removing $h$-bars until it is not possible to remove any more, and the \emph{$h$-bar-weight} of $\al$ is the number of $h$-bars removed.

We also use \abds for \bps. In this paper, we use the convention employed by Kleshchev and Livesey \cite{kl} (a different version of the abacus for \bps was introduced in \cite{yat,mfspinwt2}). We take an abacus with $h$ vertical runners labelled $0,\dots,h-1$ from left to right, with positions $0,1,2,\dots$ marked from left to right along successive rows. Given an \bp $\al$, we place a bead on the abacus at position $\al_r$ for each $r$. In particular, on runner $0$ there can be more than one bead in a given position, and we regard position $0$ as containing infinitely many beads. We call the resulting configuration the \emph{bar-\abd} for $\al$. We will use white beads when drawing bar-\abds, and we will decorate a bead on runner $0$ with $a\in\bbn\cup\{\infty\}$ to indicate that there are $a$ beads at that position.

For example, if $h=5$ and $\al=(18,12,10^2,9,7,6,2)$, we obtain the following \babd.
\[
\abacus(lmmmr,o-e\infi nonn,noono,o-e2nonn,nnnon)
\]

As with \abds, bar-\abds are useful for visualising some of the combinatorial concepts described above for \bps.
\begin{enumerate}
\item
Suppose $i\in\{1,\dots,n\}$. Then $i$-bar-addable nodes of $\al$ correspond to empty positions $b$ on runner $i+1$ or runner $h-i$ for which position $b-1$ is occupied, and $i$-bar-removable nodes correspond to occupied positions $b$ on these runners for which position $b-1$ is empty. A similar but more complicated statement holds for $i=0$.
\item
Removing an $h$-bar from $\al$ corresponds to either moving a bead from position $b\gs h$ to position $b-h$ (which must be empty if $b\nequiv0\ppmod h$), or removing beads from positions $i$ and $h-i$ for some $i\in\{1,\dots,n\}$.
\end{enumerate}

\subsection{The Fock space in type $A^{(1)}$}

Now we introduce some quantum algebra. To begin with, for any $r\in\bbn$ and an indeterminate $x$, we define the quantum integer $[r]_x=(x^r-x^{-r})/(x-x^{-1})$ and the quantum factorial $[r]^!_x=[r]_x[r-1]_x\dots[1]_x$.

Now fix an integer $m\gs2$ and an indeterminate $q$. We define $U_m$ to be the quantum group $U_{q^2}(A_{m-1}^{(1)})$ defined over $\bbc(q)$, with standard generators $\e i,\f i,\t i$ for $i\in\zmz$. We define $U^-_m$ to be the negative part of $U_m$, generated by $\lset{\f i}{i\in\zmz}$. Note in particular that we define $U_m$ with quantum parameter $q^2$ rather than $q$; this will makes the comparison with the Fock space of type $A^{(2)}$ easier.

We will be working with the level $1$ Fock space for $U_m$, which we denote $\fo$. This was introduced by Hayashi \cite{hay}, but we use the combinatorial description given by Mathas in \cite{mathas}. The Fock space $\fo$ is a $\bbc(q)$-vector space with the set $\calp$ of all partitions as a basis, which we call the \emph{standard basis}. We write $(\ ,\ )$ for the inner product on $\fo$ for which the standard basis elements are orthonormal.

Each $\la$ is a weight vector, i.e.\ a simultaneous eigenvector for the generators $\t i$. We will not need to describe the weights explicitly, but we note that two standard basis vectors $\la$ and $\mu$ have the same weight \iff  $\la$ and $\mu$ have the same $m$-core and $m$-weight.

For the purposes of this paper we only need to describe the action of the divided powers $\f i^{(r)}=\f i^r/[r]^!_{q^2}$ on the standard basis. Recall that we write $\la\arr{i:r}\mu$ if $\mu$ is obtained from~$\la$ by adding $r$ nodes of residue $i$. If this is the case, $n(\la,\mu)$ to be the sum, over all nodes $\fka$ of $\mu\setminus\la$, of the number of $i$-addable nodes of~$\mu$ to the left of $\fka$ minus the number of $i$-removable nodes of~$\la$ to the left of $\fka$. Now the action of $\f i^{(r)}$ on $\fo$ is defined by
\[
\f i^{(r)}\la=\sum_{\substack{\mu\in\calp\\\la\arr{i:r}\mu}} q^{2n(\la,\mu)}\mu.
\]
(Note in particular the factor of $2$, which arises because we take the defining quantum parameter of $U_m$ to be $q^2$, not $q$.)

We will often read the coefficient $n(\la,\mu)$ from the \abd for $\la$. To do this, we recall from above that on an $s$-bead \abd for $\la$, the addable and removable $i$-nodes can be read by looking at runners $a$ and $a-1$, where $i=a-s+n\bbz$; addable and removable nodes further to the left correspond to positions higher up these runners.

For example, suppose $m=3$ and $\la=(6,5,4,1^2)$, and take $i=1+3\bbz$ and $r=2$. Then $\la$ has three $i$-addable nodes $(2,6)$, $(4,2)$ and $(6,1)$, and one $i$-removable node $(3,4)$. So we obtain
\[
\f i^{(2)}\la=(6,5,4,2,1^2)+(6^2,4,1^3)+q^2(6^2,4,2,1).
\]
We can see this calculation from the following \abds, where $i$-addable and -removable nodes can be seen on the two leftmost runners.
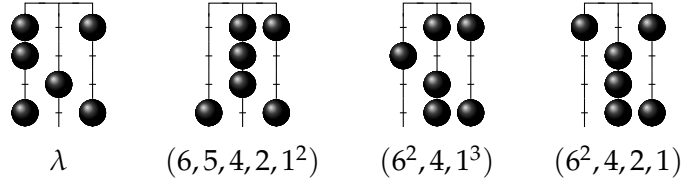
\[
\begin{array}{c@{\qquad}c@{\qquad}c@{\qquad}c}
\abacus(lmr,bnb,bnn,nbn,bnb)
&
\abacus(lmr,nbb,nbn,nbn,bnb)
&
\abacus(lmr,nbb,bnn,nbn,nbb)
&
\abacus(lmr,bnb,nbn,nbn,nbb)
\\[22pt]
\la
&
(6,5,4,2,1^2)
&
(6^2,4,1^3)
&
(6^2,4,2,1)
\end{array}
\]

Now we come to the canonical basis. Let $\calv_m$ denote the $U^-_m$-submodule generated by $\vn$ (which is in fact also the $U_m$-submodule generated by $\vn$). This is an irreducible highest-weight module for $U_m$ with highest weight $\La_0$, and admits an important involution $v\mapsto\ol v$ called the \emph{bar involution}. This is $\bbc(q+q^{-1})$-linear, and can be defined on $\calv_m$ by the properties
\[
\ol{\vn}=\vn,\qquad\ol{{\f i}v}=\f i\ol v\text{ for $i\in\zmz$}.
\]
We say that an element $v\in\calv_m$ is \emph{bar-invariant} if $\ol v=v$. The above properties imply that any $\bbc(q+q^{-1})$-linear combination of vectors of the form $\f{i_1}\dots\f{i_r}\vn$ is bar-invariant.

The bar involution allows us to define the \emph{canonical basis} of $\calv_m$, which is the main object of study in this paper. This basis is written $\lset{\cb m(\mu)}{\mu\text{ an $m$-restricted partition}}$, and the canonical basis vectors $\cb m(\mu)$ have the following properties (which are sufficient to define them uniquely).
\begin{itemize}
\item
$\cb m(\mu)$ is bar-invariant.
\item
$\cb m(\mu)$ has the form $\sum_\la \dn\la\mu\la$, where $\dn\mu\mu=1$ and $\dn\la\mu$ is a polynomial divisible by $q^2$ for $\la\neq\mu$.
\item
If $\la$ and $\mu$ are partitions and $\mu$ is $m$-restricted, then $\dn\la\mu=0$ unless $\la\dom\mu$ and $\la$ has the same $m$-core as $\mu$.
\end{itemize}
In particular, the last condition ensures that $\cb m(\mu)$ is a weight vector in $\calv_m$, so we can talk about the canonical basis of a given weight space.

The coefficients $\dn\la\mu$ are called \emph{$q$-decomposition numbers} (or in our case $q^2$-decomposition numbers), in view of Ariki's theorem that their evaluations at $q=1$ yield decomposition numbers for Iwahori--Hecke algebras at an $m$th root of unity in $\bbc$. (Of course, $\dn\la\mu$ depends on the choice of $m$, but we will always make it clear from the context which value of $m$ is intended.)

\subsection{The Fock space in type $A^{(2)}$}\label{fock2sec}

Now let $h\gs3$ be an odd integer, and set $n=\frac12(h-1)$. Define $\spin U_h$ to be the quantum group $U_q(A_{h-1}^{(2)})$. We write the standard generators for $\spin U$ as $\spe i,\spf i,\spt i$, for $i\in\{0,\dots,n\}$. We define $\spin U^-_h$ to be the $\bbc(q)$-subalgebra generated by $\spf0,\dots,\spf n$.

We let $\spin\fo$ denote the level $1$ Fock space for $\spin U$. This is a $\bbc(q)$-vector space with the set $\hstr$ of $h$-strict partitions as its standard basis. As in type $A^{(1)}$, we write $(\ ,\ )$ for the inner product with respect to which this basis is orthonormal. Each $\al$ is a weight vector, and two vectors $\al$ and $\be$ have the same weight \iff $\al$ and $\be$ have the same $h$-bar-core and $h$-bar-weight.

The action of $\spin U_h$ on $\spin\fo$ is more complicated than for $\fo$. We take the definitions from the paper \cite{lt} by Leclerc and Thibon, with additional detail from the author's paper \cite{mfspinwt2}. (Note that we use the more standard labelling of the Dynkin diagram in type $A^{(2)}_{h-1}$, so that Leclerc and Thibon's generators $\f0,\dots,\f n$ are our $\spf n,\dots,\spf0$).

Given $i\in\{0,\dots,n\}$ and $r\in\bbn$, we define $\spf i^{(r)}=\spf i^r/[r]^!_{q_i}$, where
\[
q_i=
\begin{cases}
q&\text{if }i=0\\
q^2&\text{if }0<i<n\\
q^4&\text{if }i=n.
\end{cases}
\]
Recall that we write $\al\aarr{i:r}\be$ if $\be$ is obtained from $\al$ by adding $r$ nodes of bar-residue $i$. If this is the case, we define $\spin n(\al,\be)$ to be the sum, over all nodes $\fka$ of $\be\setminus\al$, of the number of $i$-bar-addable nodes of $\be$ to the left of $\fka$ minus the number of $i$-bar-removable nodes of $\al$ to the left of $\fka$. Further, if $i=0$, let $M$ be the set of integers $m\gs1$ such that column $mh+1$ contains a node of $\be\setminus\al$ but column $mh$ does not. For each $m\in M$, let~$b_m$ be the number of times $mh$ occurs as a part of $\al$, and set $N=\prod_{m\in M}(1-(-q^2)^{b_m})$. (If $i\neq0$, then set $N=1$.)

Now the action of $\spf i$ is given by
\[
\spf i^{(r)}\al=\sum_{\substack{\be\in\hstr\\\al\aarr{i:r}\be}}Nq_i^{\spin n(\al,\be)}\be.
\]
As in type $A^{(1)}$, we will often read the coefficient $\spin n(\al,\be)$ from the \babds for $\al$ and $\be$; here the situation is slightly more complicated when $i<n$, since there are more than two runners to consider.

Now we can define the canonical basis; this is done in essentially the same way as in type $A^{(1)}$. Let $\spin\calv_h$ denote the $\spin U^-_h$-submodule generated by $\vn$ (which is also the $\spin U_h$-submodule generated by $\vn$). This is an irreducible module with highest weight $\La_0$. The \emph{bar involution} on $\spin\calv_h$ is the $\bbc(q+q^{-1})$-linear map determined by
\[
\ol{\vn}=\vn,\qquad\ol{{\spf i}v}=\spf i\ol v\text{ for $i\in\{0,\dots,n\}$}.
\]

The canonical basis for $\spin\calv_h$ is written
\[
\lset{\scb h(\be)}{\be\text{ a restricted $h$-strict partition}},
\]
and the canonical basis vectors $\scb h(\be)$ are defined by the following properties.
\begin{itemize}
\item
$\scb h(\be)$ is bar-invariant.
\item
$\scb h(\be)$ has the form $\sum_\al\sdn\al\be\al$, where $\sdn\be\be=1$ and $\sdn\al\be$ is a polynomial divisible by $q$ for $\al\neq\be$.
\item
$\sdn\al\be=0$ unless $\al\dom\be$ and $\al$ has the same $h$-bar-core as $\be$.
\end{itemize}
Our main focus in this paper is comparing canonical bases in types $A^{(1)}$ and $A^{(2)}$.

\section{Symmetric functions}\label{symfnsec}

In this section we recall some basic theory of symmetric functions, and prove an apparently new result which we will use later.

\subsection{Background on symmetric functions}

We take an indeterminate $t$ over $\bbq$, and a countably infinite set $X$ of commuting indeterminates. We let $\La$ be the ring of \emph{symmetric functions}: power series in the elements of $X$ with bounded degree with coefficients in $\bbq(t)$, which are invariant under permutations of $X$.

$\La$ has several important bases. Two of the most important are the bases of \emph{Schur functions} $s_\la$ and \emph{Hall--Littlewood functions} $P_\la$. We refer to Macdonald's book \cite{macd} for definitions of these, as well as a detailed introduction to symmetric functions.

We will need to use the standard coproduct $\Delta:\La\to\La\otimes\La$. To define this, first note that for any countably infinite set $Y$ and any $f\in\La$ we can define $f(Y)$ simply by replacing the elements of $X$ with their images under some chosen bijection from $X$ to $Y$. If we partition the set of variables $X$ into two disjoint infinite sets $Y\sqcup Z$, then each symmetric function $f$ is symmetric in the elements of $Y$ and the elements of $Z$, so we can write $f=\sum_{i\in I}g_i(Y)h_i(Z)$ for some finite indexing set $I$ and $f_i,g_i\in\La$; we then define $\Delta(f)=\sum_ig_i\otimes h_i$.

Our aim is to study Pieri-type rules for Hall--Littlewood functions, and for this we need some notation. Suppose $\la,\mu\in\calp$. Recall that we write $\la\car\mu$ if $\la\subseteq\mu$ and $\mu\setminus\la$ consists of $r$ nodes lying in different columns. (In the literature on symmetric functions, this is often expressed as saying that $\mu\setminus\la$ is a \emph{horizontal strip} of length $r$.) If this is the case, then we define
\[
\hs\la\mu=\sum_c(1-t^{\mu'_c-\mu'_{c+1}}),
\]
summing over all $c\gs1$ such that column $c$ contains a node of $\skw\mu\la$ but column $c+1$ does not. The polynomials $\hs\la\mu$ are used to write down a ``Pieri rule'' for Hall--Littlewood functions. To define this, we need the second type of Hall--Littlewood functions $Q_\la$ introduced by Macdonald. For any $n\in\bbn$ define
\begin{align*}
\phi_n&=\prod_{i=1}^n(1-t^i).
\\
\intertext{Now for any partition $\la$ define}
b_\la&=\prod_c\phi_{\la'_c-\la'_{c+1}}.
\\
\intertext{Then $Q_\la$ is defined by}
Q_\la&=b_\la P_\la.
\end{align*}
Now we can state Macdonald's Pieri rule for Hall--Littlewood functions.

\begin{thm}[\xcite{macd}{III.5.7}]\label{hlpieri}
If $\la\in\calp$ and $r\gs0$ then
\[
P_\la Q_{(r)}=\sum_{\substack{\mu\in\calp\\\la\car\mu}}\hs\la\mu P_\mu.
\]
\end{thm}

\subsection{A new Pieri rule}\label{newpierisec}

Our objective in this section is to prove a kind of dual rule to \cref{hlpieri}, which surprisingly seems not to be in the literature.

Suppose $f\in\La$. Because the Schur functions $s_\la$ form a basis for $\La$ we can uniquely write $\Delta(f)=\sum_\la\ddd\la(f)\otimes s_\la$ for some symmetric functions $\ddd\la(f)$. This defines a function $\ddd\la:\La\to\La$ for each partition $\la$. We are particularly interested in the functions $\ddd{(r)}$ for $r\in\bbn$. One form of the classical Pieri rule is that $\ddd{(r)}s_\mu=\sum s_\la$, summing over all $\la$ such that $\la\car\mu$. To give the corresponding rule for $\ddd{(r)}P_\mu$, we need some more notation. If $\la,\mu$ are partitions with $\la\car\mu$, define
\[
\hsb\la\mu=\sum_c(1-t^{\la'_c-\la'_{c+1}}),
\]
summing over all $c\gs1$ such that column $c+1$ contains a node of $\skw\la\mu$ but column $c$ does not.

The main result of this section is the following, which generalises part of \cite[Proposition~3.6]{mfspin2alt}.

\begin{propn}\label{dualpieri}
Suppose $\mu\in\calp$ and $r\in\bbn$. Then
\[
\ddd{(r)}P_\mu=\sum_{\substack{\la\in\calp\\\la\car\mu}}\hsb\la\mu P_\la.
\]
\end{propn}

We will show how to derive \cref{dualpieri} from \cref{hlpieri} using the self-duality of $\La$ explained by Konvalinka and Lauve in \cite[\S2]{konlau}. The \emph{Hall inner product} on $\La$ is the bilinear function defined by $\lan P_\mu,Q_\la\ran=\delta_{\mu\la}$. By \cite[Lemma 11]{konlau}, $\La$ is self-dual with respect to $\lan\,,\,\ran$, which means that (defining $\lan\,,\,\ran$ on $\La\otimes\La$ as well in the obvious way)
\[
\lan\Delta(f),g\otimes h\ran=\lan f,gh\ran
\]
for any $f,g,h\in\La$.

The final ingredient we need is the following \lcnamecref{bstrip}, which follows immediately from the definitions.

\begin{lemma}\label{bstrip}
Suppose $\la,\mu\in\calp$ and $\la\car\mu$. Then
\[
\hsb\la\mu=\frac{b_\la\hs\la\mu}{b_\mu}.
\]
\end{lemma}

To complete the proof of \cref{dualpieri}, we need to consider the transition matrix between the bases $(s_\la)_{\la\in\calp}$ and $(P_\la)_{\la\in\calp}$. This matrix is denoted $K(t)$, and its entries $K_{\la\mu}(t)$ are called \emph{Kostka polynomials}. Specifically, we write
\[
s_\la=\sum_\mu K_{\la\mu}(t)P_\mu,\qquad P_\la=\sum_\mu K^{-1}_{\la\mu}(t)s_\mu.
\]

We will need the following properties of these polynomials.

\begin{lemma}\label{kostkaprops}
Suppose $\la,\mu\in\calp$. Then $K_{\la\mu}(t)$ and $K^{-1}_{\la\mu}(t)$ are polynomials in $t$ which are zero unless $\la\dom\mu$. Furthermore, $K_{\la\la}(t)=K^{-1}_{\la\la}(t)=1$, while $K_{\la\mu}(t)$ and $K^{-1}_{\la\mu}(t)$ are divisible by $t$ when $\la\neq\mu$.
\end{lemma}

\begin{pf}
The given properties for $K_{\la\mu}(t)$ follow from \cite[III.2.3 \& III.2.6]{macd}. Now the properties for $K^{-1}_{\la\mu}(t)$ follow by inverting the matrix $K(t)$.
\end{pf}

\begin{pf}[Proof of \cref{dualpieri}]
Let $\lset{S_\la}{\la\in\calp}$ be the basis dual to the basis of Schur functions with respect to $\lan\ ,\ \ran$. Because $\lset{Q_\la}{\la\in\calp}$ and $\lset{s_\la}{\la\in\calp}$ are both bases for $\La$, there are coefficients $a_{\mu\pi\nu}$ defined for all $\mu,\pi,\nu\in\calp$ such that
\[
\Delta(P_\mu)=\sum_{\pi,\nu\in\calp}a_{\mu\pi\nu}Q_\pi\otimes s_\nu.
\]
Then
\[
a_{\mu\pi\nu}=\lan\Delta(P_\mu),P_\pi\otimes S_\nu\ran.
\]
By \cref{kostkaprops} we can write $s_\la$ as $P_\la$ plus a linear combination of the $P_\mu$ with $\la\doms\mu$, for each $\la$. As a consequence, $S_{(r)}$ equals $Q_{(r)}$. So
\begin{align*}
\ddd {(r)}P_\mu
&=\sum_{\pi\in\calp} a_{\mu\pi(r)}Q_\pi\\
&=\sum_{\pi\in\calp}\lan\Delta(P_\mu),P_\pi\otimes Q_{(r)}\ran Q_\pi\\
&=\sum_{\pi\in\calp}\lan P_\mu,P_\pi Q_{(r)}\ran b_\pi P_\pi\\
&=\sum_{\substack{\pi,\la\in\calp\\\pi\car\la}}\frac{\hs\pi\la}{b_\mu}\lan Q_\mu,P_\la\ran b_\pi P_\pi\tag*{by \cref{hlpieri}}\\
&=\sum_{\substack{\pi\in\calp\\\pi\car\mu}}\frac{\hs\pi\mu b_\pi}{b_\mu}P_\pi\\
&=\sum_{\substack{\pi\in\calp\\\pi\car\mu}}\hsb\pi\mu P_\pi.\tag*{by \cref{bstrip}.\qedhere}
\end{align*}
\end{pf}

\section{Comparing canonical bases in types $A^{(1)}_{n-1}$ and $A^{(2)}_{h-1}$}\label{firstcomparesec}

\begin{framing}{black}
For the rest of the paper we fix a natural number $l$. We fix an odd integer $h\gs3$, and let $n=\frac12(h-1)$ and $m=\frac12(h+1)$.
\end{framing}

Our aim in this paper is to compare the canonical bases in $\calv_m$ and $\spin\calv_h$. But as an intermediate step in this section, we compare $\calv_n$ with $\spin\calv_h$. So \emph{for this section, we assume that $h\gs5$}.

We begin by setting up some combinatorics underlying our comparison between canonical bases. Say that an \bp $\al$ is \emph{\nice} if the residue modulo $h$ of every non-zero part of $\al$ lies in $\{1,\dots,n\}$. Let $\pnice l$ denote the set of \nice \bps of length $l$, and let $\pl l$ denote the set of all partitions of length at most $l$. Proving the next \lcnamecref{nicebijec} is a routine exercise.

\begin{propn}\label{nicebijec}
Fix $l\gs1$. There is a bijection $\phi:\pnice l\to\pl l$ given by
\begin{align*}
(\phi\al)_r&=\al_r-(n+1)\inp{\frac{\al_r}h}-l+r-1
\\
\intertext{for $1\ls r\ls l$. The inverse of $\phi$ is given by}
(\phi^{-1}\la)_r&=\la_r+l-r+1+(n+1)\inp{\frac{\la_r+l-r}n}.
\end{align*}
If $\be\in\pnice l$, then $\be$ is restricted \iff $\phi\be$ is $n$-restricted.
\end{propn}

The bijection $\phi$ is easily realised on the abacus. Given $\al\in\pnice l$, the \babd for $\al$ has $l$ beads on runners $1,\dots,n$, no beads on runners $n+1,\dots,h-1$, and no beads on runner $0$ except in position $0$. Deleting runners $0$ and $n+1,\dots,h-1$ yields the \abd for $\phi\al$.

For example, taking $h=5$ and $\al=(16,12,11,6,2)$, we obtain $\phi\al=(2^3,1^2)$.
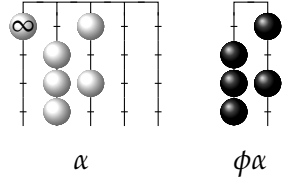
\[
\begin{array}{c@{\qquad}c}
\abacus(lmmmr,o-e\infi nonn,nonnn,noonn,nonnn)
&
\abacus(lr,nb,bn,bb,bn)
\\[22pt]
\al&\phi\al
\end{array}
\]

We want to compare canonical basis elements labelled by partitions in $\pnice l$ and $\calp_{\ls l}$. So let
\[
\fo^{\ls l}=\lspan{\la}{\la\in\pl l},\qquad\fonicel=\lspan{\al}{\al\in\pnice l},
\]
where $\lan X\ran$ denotes the $\bbc(q)$-span of a set $X$. Now we have a bijective linear map $\Phi:\fonicel\to\fo^{\ls l}$ defined by mapping $\al\mapsto\phi\al$ and extending linearly.

Observe that if $\mu\in\pl l$ is $n$-restricted, then $\cb n(\mu)\in\fo^{\ls l}$, because if $\la$ appears in $\cb n(\mu)$ with non-zero coefficient then $\la\dom\mu$, so that $l(\la)\ls l$. A similar statement is true in $\fonicel$, but requires the following lemma.

\begin{lemma}\label{nicecore}
Suppose $\be\in\pnice l$. If $\al$ is an \bp with the same $h$-bar-core as $\be$ and $\al\dom\be$, then $\al\in\pnice l$.
\end{lemma}

\begin{pf}
First we observe that if we remove an $h$-bar from $\be$, the resulting \bp will be \nice with length $l$: since $\be$ is \nice, it does not have two parts which sum to $h$, and so removing an $h$-bar from $\be$ entails reducing by $h$ a part which is larger than $h$. This does not change the length of the partition, or the set of residues modulo $h$ of the non-zero parts. So the resulting partition is \nice with length $l$. Applying this repeatedly, we find that the $h$-bar-core $\ga$ of $\be$ is \nice with length $l$.

Now consider $\al$. Because $\al\dom\be$, the length of $\al$ is at most $l$. On the other hand, by assumption the $h$-bar-core of $\al$ is $\ga$ which has length $l$, and therefore $l(\al)=l$. We can obtain $\al$ from $\ga$ by repeatedly adding $h$-bars. At each stage we do not increase the length of the partition, so addition of an $h$-bar must consist of increasing some positive part by $h$. This does not affect the set of residues modulo $h$ of the non-zero parts, so $\al$ is \nice.
\end{pf}

As a consequence, if $\be\in\pnice l$ is restricted, then $\scb h(\be)\in\fonicel$. Our first main result is the following.

\begin{thm}\label{samedec}
Suppose $\beta\in\pnice l$ is restricted. Then
\[
\cb n(\phi\beta)=\Phi(\scb h(\beta)).
\]
\end{thm}

\begin{eg}
Take $h=5$, so that $n=2$. Let $\be=(16,12,11,7,6,1)$. Then $\phi\be=(1^5)$. Writing $\f i=\f{i+2\bbz}$ for $i=0,1$, we can calculate
\[
\big(\f0\f1\f0\f1\f0-\f0^{(2)}\f1^{(2)}\f0\big)\vn=(1^5)+q^2(3,1^2)+q^4(5),
\]
so that
\[
\cb2(1^5)=(1^5)+q^2(3,1^2)+q^4(5).
\]
On the other hand, $(12,11,7,6,2,1)$ is a $5$-bar-core, so $(12,11,7,6,2,1)=\scb5(12,11,7,6,2,1)$ is bar-invariant. If we write $\spin g=\spf0^{(2)}\spf1\spf2$ and $\spin g^{(2)}=\spin g^2/(q^2+q^{-2})$, then
\[
\big(\spin g\spf1\spin g\spf1\spin g-\spin g^{(2)}\spf1^{(2)}\spin g)(12,11,7,6,2,1)=
(16,12,11,7,6,1)+q^2(21,12,11,6,2,1)+q^4(26,11,7,6,2,1)
\]
is also bar-invariant, so equals $\scb5(\be)$. So $\Phi(\scb5(\be))=\cb2(\phi\beta)$.
\end{eg}

The above example suggests how \cref{samedec} is proved: we define linear operators $g_0,\dots,g_{n-1}$ on $\spin\fo$ which correspond to the action of the generators $\f i$ on $\fo$.

For the next two results we take $i\in\{1,\dots,n-1\}$, and let $\hi=i-l+n\bbz$.

\begin{lemma}\label{simpleaddres}
Suppose $\al\in\pnice l$ and $1\ls r\ls l$. Then $\al$ has an $i$-bar-addable node $\fkb$ in row $r$ \iff $\phi\al$ has an $\hi$-addable node $\fkc$ in row $r$. If these nodes exist, then $\al\cup\fkb$ lies in $\pnice l$, and $\phi(\al\cup\fkb)=\phi\al\cup\fkc$.
\end{lemma}

\begin{pf}
Since $\al$ is \nice, $\al_r$ cannot be congruent to $h-i-1$ modulo $h$, so $\al$ has an $\hi$-bar-addable node $\fkb$ in row $r$ \iff $\al_r\equiv i\ppmod h$ and either $r=1$ or $\la_{r-1}-\la_r\gs2$. Applying $\phi$, we see that this is equivalent to the condition that $\phi\al$ has an $\hi$-addable node $\fkc$ in row $r$.

If the nodes $\fkb$ and $\fkc$ do exist, then $(\al\cup\fkb)_r=\al_r+1\equiv\hi+1\ppmod h$, so that $\al\cup\fkb$ is \nice. Furthermore, $\phi(\al\cup\fkb)_r=(\phi\al)_r+1$, so $\phi(\al\cup\fkb)=\phi\al\cup\fkc$.
\end{pf}

A corresponding result for $\hi$-bar-removable nodes can be proved in exactly the same way. As a consequence, we obtain the following.

\begin{cory}\label{simpleg}
If $\al$ is a \nice \bp, then $\spf i\al\in\fonicel$, and
\[
\Phi(\spf i\al)=\f\hi\phi\al.
\]
\end{cory}

\begin{pf}
Since $1\ls i\ls n-1$, adding an $i$-bar-addable node to a \nice \bp yields a \nice \bp; so the first statement holds.

Write $\la=\phi\al$. Then $\f\hi\la$ is a linear combination of basis elements $\mu$ obtained by adding a $\hi$-addable node to $\la$. Take such a $\mu$, and suppose the added node lies in row $r$. Then $r\ls l$, since the residue of the addable node of $\la$ in row $l+1$ (if there is one) is $-l+n\bbz\neq\hi$. So $\mu=\phi\be$, where $\be\in\pnice l$ is obtained from $\al$ by adding a node in row $r$. The coefficient of $\mu$ in $\f\hi\la$ is determined by the $\hi$-addable and $\hi$-removable nodes of $\la$ below row $r$. These nodes all lie in row $l$ or higher, so by \cref{simpleaddres} (and its analogue for removable nodes) they correspond to the $i$-bar-addable and $i$-bar-removable nodes of $\al$ in rows $r+1,\dots,l$. Furthermore, $\al$ does not have an $i$-bar-addable node below row $l$ (because $i\neq0$). So (comparing the action of $\f\hi$ on $\fo$ with the action of $\spf i$ on $\spin\fo$) the coefficients $\ipx{\spf i\al}{\be}$ and $\ip{\f\hi\la}{\mu}$ agree.
\end{pf}

Now we consider the remaining residue not included in the last two results: let $\hi=-l+n\bbz$, and let $\spin g=\spf0^{(2)}\spf 1\dots\spf n$.

\begin{propn}\label{hardg}
Suppose $\al\in\pnice l$, and $\be$ is an \bp with length $l$. Then \tfae.
\begin{enumerate}
\item
$\ip{\spin g\al}{\be}\neq0$.
\item
$\be$ is \nice and $\ip{\f\hi\phi\al}{\phi\be}\neq0$.
\end{enumerate}
Furthermore, if these two conditions hold, then $\ip{\spin g\al}{\be}=\ip{\f\hi\phi\al}{\phi\be}$.
\end{propn}

\begin{pf}
Suppose $\be$ appears with non-zero coefficient in $\spin g\al$. Then there is a sequence of \bps $\al(1),\dots,\al(n)$ such that
\[
\al\aarr n\al(1)\aarr{n-1}\cdots\aarr1\al(n)\aarr{0:2}\be.
\]
In particular, $\be$ has two nodes of bar-residue $0$ that are not contained in $\al$. Let $r$ be a row containing at least one of these nodes; then the node equals $(r,hc)$ or $(r,hc+1)$ for some $c\gs1$. This means that $\al_r\ls hc$, but because $\al$ is \nice, this actually gives $\al_r\ls hc-n-1$, so that the nodes
\[
(r,hc- n),\dots,(r,hc)
\]
all belong to $\be$ but not to $\al$. The same applies for the other node of bar-residue $0$ that lies in $\be$ but not in~$\al$, so in fact these nodes must lie in the same row. Hence there is some $r\ls l$ such that $\al_r\equiv n\ppmod h$ and
\[
\be_s=
\begin{cases}
\al_s+n+2&\text{if }s=r
\\
\al_s&\text{otherwise}.
\end{cases}
\]
In particular, $\be$ is \nice. Furthermore, applying the definition of $\phi$ we see that $\phi\be=\phi\al\cup\fkn$, where $\fkn$ is an $\hi$-addable node of $\phi\al$. So (1)$\Rightarrow$(2), and the converse is similar (but easier).

It remains to compare the coefficient of $\be$ in $\spin g\al$ with the coefficient of $\phi\be$ in $\f\hi\phi\al$. To compute the coefficient of $\be$ in $\spin g\al$, we consider bar-addable and bar-removable nodes below row $r$. The coefficient $\ip{\f\hi\phi\la}{\phi\be}$ includes a factor $q^2$ for each bar-addable node of $\al$ below row $r$, and a factor $q^{-2}$ for each bar-removable node, except for $n$-bar-addable nodes which give a factor of $q^4$ (note that $\al$ has no $n$-bar-removable nodes because it is \nice). The fact that $\al$ is \nice also means that for each $s>r$, the number of bar-removable nodes in row $s$ equals the number of bar-addable nodes in row $s+1$, except when $\al_s\equiv1\ppmod h$, in which case there is one more bar-removable node in row $s$. Furthermore, there is an $n$-bar-addable node in row $s$ \iff $\al_s\equiv n\ppmod h$. So if we let $a_1$ be the number of $s>r$ such that $\al_s\equiv1\ppmod h$, and define $a_n$ similarly, then
\[
\ip{\spin g\al}{\be}=q^{2(a_n-a_1+1)}.
\]
(The extra $1$ arises from the bar-addable node in row $r+1$.)

But now observe that $\al_s\equiv n\ppmod h$ \iff $(\phi\al)_s-s\in\hi-1$, while $\al_s\equiv1\ppmod h$ \iff $(\phi\al)_s-s\in\hi$. So $a_n-a_1+1$ equals the number of $i$-addable nodes of $\phi\la$ below row $r$ minus the number of $\hi$-removable nodes (the extra $1$ arises from the addable node in row $l+1$). So $\ip{\spin g\la}{\be}=\ip{\f\hi\phi\al}{\phi\be}$.
\end{pf}

This gives an analogue of \cref{simpleg}. We continue to write $\hi=-l+n\bbz$.

\begin{cory}\label{hardq}
Suppose $\al$ is a \nice \bp, and that $\f\hi\phi\al\in\fo^{\ls l}$. Then $\spin g\al\in\fonicel$, and
\[
\Phi(\spin g\al)=\f\hi\phi\al.
\]
\end{cory}

\begin{pf}
From the fact that (1)$\Rightarrow$(2) in \cref{hardg} we can write
\[
\spin g\al=\sum_{\be\in\pnice l}c_\be\be
\]
for some coefficients $c_\be\in\bbc(q)$, so that $\spin g\al\in\fonicel$. Then the last statement in \cref{hardg} tells us that for each $\mu\in\pl l$ the coefficient of $\mu$ in $\f\hi\phi\al$ is $c_{\phi^{-1}\mu}$. By assumption no $\mu$ with $\mu\notin\pl l$ occurs in $\f\hi\phi\al$, and so
\[
\f i\phi\al=\sum_{\mu\in\pl l}c_{\phi^{-1}\mu}\mu=\sum_{\be\in\pnice l}c_\be\Phi(\be)=\Phi(\spin g\al).\qedhere
\]
\end{pf}

Now we can prove our first main result.

\begin{pf}[Proof of \cref{samedec}]
Take a restricted partition $\be\in\pnice l$, and let $\mu=\phi\be$. The LLT algorithm \cite[Section~6.2]{llt} shows that $\cb n(\mu)$ can be written as a linear combination $\sum_\nu a_\nu A(\nu)$, where:
\begin{enumerate}
\item\label{nudommu}
the sum is over a set of partitions $\nu$ with $\nu\dom\mu$;
\item
each coefficient $a_\nu$ lies in $\bbc(q+q^{-1})$;
\item
each $A(\nu)$ has the form $\f{i_1}\dots \f{i_r}\vn$ for some $i_1,\dots,i_r\in\znz$;
\item\label{xidomnu}
each $A(\nu)$ is a linear combination of partitions $\xi$ with $\xi\dom\nu$.
\end{enumerate}
In particular, conditions (\ref{nudommu}) and (\ref{xidomnu}) imply that each $A(\nu)$ lies in $\fo^{\ls l}$. So if we take such a $\nu$ and write
\[
A(\nu)=\f{i_1}\dots \f{i_r}\vn=\sum_{\xi\in\pl l}c_{\nu\xi}\xi
\]
with each $c_{\nu\xi}\in\bbc(q)$, then by \cref{simpleg,hardq} we can find $j_1,\dots,j_s\in\{0,\dots,n\}$ and a coefficient $b_\nu\in\bbc(q+q^{-1})$ such that
\[
\Phi^{-1}(A(\nu))=b_\nu\spf{j_1}\dots\spf{j_s}\phi^{-1}\vn=\sum_{\xi\in\pl l}c_{\nu\xi}\phi^{-1}\xi.
\]
(The coefficient $b_\nu$ arises because of the divided power occurring in the definition of the operator $\spin g$.) Hence the vector
\[
\Phi^{-1}(\cb n(\mu))=\sum_\nu a_\nu\Phi^{-1}(A(\nu))
\]
can be written as a linear combination of vectors of the form $\spf{j_1}\dots\spf{j_s}\phi^{-1}\vn$, with coefficients in $\bbc(q+q^{-1})$. Now $\phi^{-1}\vn$ is an $h$-bar-core: its parts are simply the smallest $l$ positive integers whose residues modulo $h$ lie in $\{1,\dots,n\}$. So $\phi^{-1}\vn=\scb h(\phi^{-1}\vn)$ is bar-invariant, and so $\Phi^{-1}(\cb n(\mu))$ is bar-invariant. Since $\Phi^{-1}(\cb n(\mu))=\sum_\la \dn\la\mu\phi^{-1}\la$, with $\dn\mu\mu=1$ and all other $\dn\la\mu$ divisible by $q$, the uniqueness of the canonical basis means that $\Phi^{-1}(\cb n(\mu))=\scb h(\be)$.
\end{pf}

\section{Comparing types $A^{(1)}_{n-1}$ and $A^{(1)}_{m-1}$: runner addition}\label{runaddsec}

In this section we invoke a theorem from the author's paper \cite{mfrunrem} comparing canonical bases in types $A^{(1)}_{n-1}$ and $A^{(1)}_{m-1}$, and use it to express \cref{samedec} in terms of the canonical basis for $\calv_m$; this will also allow us to include the case $h=3$.

Let $\pml$ denote the set of partitions $\la$ with length at most $l$ such that $\la_r+l+1-r\nequiv0\ppmod l$ for $1\ls r\ls l$. Another way of saying this is that $\la\in\pml$ if $\la$ can be displayed on an $m$-runner abacus with $l+1$ beads such that runner $0$ contains only a bead in position $0$. We define a bijection $\pl l\to\pml$ which we denote $\la\mapsto\la^+$, by setting
\[
\la^+_r=\la_r+\inp{\frac{\la_r+l-r}n}
\]
for $r=1,\dots,l$. This bijection is easily visualised on the abacus: given $\la\in\pl l$, we take the $n$-runner \abd for $\la$ with $l$ beads, add a runner at the left with a bead in the top position only, and let $\la^+$ be the resulting partition. It is easy to see that $\la^+$ is $m$-restricted \iff $\la$ is $n$-restricted. Extending the map $\la\mapsto\la^+$ linearly, we obtain an injective linear map $\fo^{\ls l}\to\fo$ written $v\mapsto v^+$.

Now we can give the main result from \cite{mfrunrem}. Some translation of notation is required, since in \cite{mfrunrem} the opposite convention for $\fo$ is used, in which the canonical basis elements $\cb n(\mu)$ are indexed by $n$-regular partitions, and a ``full'' runner is added rather than an empty one. But replacing all partitions with their conjugates yields the following result.

\begin{thm}[\xcite{mfrunrem}{Theorem 3.1}]\label{runadd}
Suppose $n\gs2$ and $\la\in\pl l$ is $n$-restricted. Then $\cb m(\la^+)=\cb n(\la)^+$.
\end{thm}

We combine this with the results of the last section to give a result comparing canonical bases in types $A^{(1)}_{m-1}$ and $A^{(2)}_{h-1}$. Recalling the bijection $\phi:\pnice l\to\pl l$ from \cref{firstcomparesec}, let $\psi:\pnice l\to\pml$ be the bijection defined by $\psi\la=(\phi\la)^+$. Explicitly, $\psi$ is given by
\[
(\psi\al)_r=\al_r-(m-1)\inp{\frac{\al_r}h}-l+r-1
\]for $1\ls r\ls l$.

We define an injective $\bbc(q)$-linear map $\Psi:\fonicel\to\fo$ by $\Psi(\be)=\psi\be$. Now we obtain the following result.

\begin{thm}\label{firstmain}
Suppose $\be\in\pnice l$ is restricted. Then
\[
\cb m(\psi\be)=\Psi(\scb h(\be)).
\]
\end{thm}

\begin{pf}
If $h\gs5$ then this follows immediately from \cref{samedec,runadd}. In the case $h=3$, the assumption that $\be$ is restricted means that $\be$ is the partition $(3l-2,3l-5,\dots,1)$, which is a $3$-bar-core, so that $\scb3(\be)=\be$. On the other hand, $\psi\be$ is the partition $(l,l-1,\dots,1)$ which is a $2$-core, so that $\cb2(\psi\be)=\psi\be$. So the result follows in this case too.
\end{pf}

\section{Comparing types $A^{(1)}_{m-1}$ and $A^{(2)}_{h-1}$: separated partitions}\label{nextbitsec}

Having established \cref{firstmain}, we will work with $U_m$ rather than $U_n$ from now on. With this in mind, if $a$ is an integer then we may write $\f a$ to mean $\f{a+m\bbz}$.

Our next aim is to extend \cref{firstmain} to include $h$-strict partitions with positive parts divisible by $h$. To do this, we need to restrict attention to what we call \emph{separated} partitions. In this section we define separated partitions and prove two preparatory results giving similar calculations in $\fo$ and $\spin\fo$.

\subsection{Separated partitions}\label{sepsec}

We keep $l\in\bbn$ fixed. Given a partition $\la\in\pml$ and any partition $\pi$, we define a partition $\jo\la\pi$ as follows: we take the \abd for $\la$ with $m$ runners and $am+l+1$ beads for sufficiently large $a$, and then move the $r$th lowest bead on runner $0$ down $\pi_r$ positions for each $r$. Another way to express this is as follows: we write $\la=\mu^+$ for $\mu\in\pl l$; then we join together the $1$-runner \abd for $\pi$ with $a+1$ beads and $n$-runner \abd for $\mu$ with $an+l$ beads. It is easy to see that $\jo\la\pi$ is independent of the choice of $a$.

Given such an \abd for $\jo\la\pi$, let $b$ be the position of the last bead on runner $0$ (that is, $b=(\pi_1+a)m$), and let $f$ be the first empty position not on runner $0$. Given $k\in\bbn_0$, we say that $\jo\la\pi$ is \emph{$k$-\sep} if $f-b>km$. For $k=0$, we will just say \emph{\sep} rather than $0$-\sep. Whether $\jo\la\pi$ is $k$-\sep is independent of the choice of $a$ (in fact, it only depends on $l$, $k$, $\pi_1$ and $\la'_1$).

\begin{eg}
Take $h=7$, so that $m=4$, and let $l=10$. If we take $\la=(7,6^2,4,1^3)$ and $\pi=(1^2)$, then $\jo\la\pi=(7,6^2,4,1^{11})$, as we see from the following \abds.
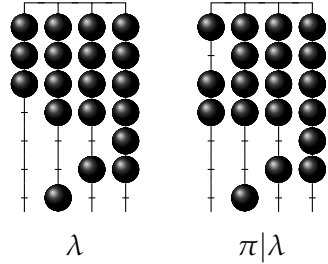
\[
\begin{array}{c@{\qquad}c}
\abacus(lmmr,bbbb,bbbb,bbbb,nbbb,nnnb,nnbb,nbnn)
&
\abacus(lmmr,bbbb,nbbb,bbbb,bbbb,nnnb,nnbb,nbnn)
\\[38pt]
\la
&
\jo\la\pi
\end{array}
\]
We can see that $\jo\la\pi$ is $1$-\sep but not $2$-\sep.
\end{eg}

Now for each $k\in\bbn$ we define an operator $F_k$ on the Fock space by
\[
F_k=\f{-l}^{(k)}\f{1-l}^{(k)}\dots \f{m-1-l}^{(k)}.
\]

The first result we want to prove describes $F_k(\jo\la\pi)$ when $\jo\la\pi$ is $k$-\sep. Recall that given two partitions $\pi,\rho$, we write $\pi\car\rho$ if $\rho$ can be obtained from $\pi$ by adding $r$ nodes in distinct columns.

\begin{propn}\label{addrun1simple}
Suppose $\la\in\pml$, $\pi\in\calp$ and $k,u\in\bbn_0$, and that $\jo\la\pi$ is $(k+u)$-\sep. If $\nu\in\calp$, then $\ip{F_k(\jo\la\pi)}{\nu}\neq0$ \iff $\nu$ has the form $\jo\mu\rho$, where
\begin{itemize}
\item
$\mu\in\pml$ and $\rho\in\calp$,
\item
$\pi\car\rho$ for some $0\ls r\ls k$, and
\item
$\ip{F_{k-r}\la}{\mu}\neq0$.
\end{itemize}
Furthermore, if these conditions hold, then $\jo\mu\rho$ is $u$-\sep.
\end{propn}

\begin{pf}
In this proof we use an abacus with $m$ runners and $l+am+1$ beads for fixed large $a$. For any partition $\nu$ we write $\ab\nu$ for the $m$-runner \abd for $\nu$ with $l+am+1$ beads. 
Given a finite set $B$ of positions on the abacus, we use the phrase \emph{making bead moves from $B$} to mean moving a bead from position $b$ to position $b+1$ for each $b\in B$.

Suppose $\nu$ appears in $F_k(\jo\la\pi)$ with non-zero coefficient. From the description of the action of the $\f i$ on $\fo_m$ (interpreted in terms of the abacus) there are sets $B_0,\dots,B_{m-1}\subset\bbn_0$ with $\card{B_g}=k$ and $B_g\subset g+m\bbz$ for each $g$, such that $\ab\nu$ is obtained from $\ab{\jo\la\pi}$ by making bead moves from $B_{m-1},B_{m-2},\dots,B_0$ in turn. We write $B=B_0\cup\dots\cup B_{m-1}$.

Let $xm$ be the position of the last bead on runner $0$ in $\ab{\jo\la\pi}$; that is, $x=\pi_1+a$. The assumption that $\jo\la\pi$ is $(k+u)$-\sep means that every position before position $(x+k+u)m$ not on runner $0$ is occupied in $\ab{\jo\la\pi}$. We make two observations about the relationships between the sets $B_0,\dots,B_{m-1}$.
\begin{enumerate}[label=(\roman*)]
\item
If $b\in B_0$ with $b>xm$, then $b-1\in B_{m-1}$ (otherwise it is not possible to make a bead move from $b$ at the final step).
\item\label{underx}
If $b\in B_g$ with $0\ls g<m-1$ and $b<(x+k+u)m$, then $b+1\in B_{g+1}$ (so that there is a space to move a bead into from position $b$ when making bead moves from $B_g$). Applying this repeatedly, we deduce that if $b\in B_0$ with $b<(x+k+u)m$, then $b+m-1\in B_{m-1}$.
\end{enumerate}
Since $\card{B_0}=k$, we can find $y\in\{x,x+1,\dots,x+k\}$ such that $ym\notin B_0$. Now the two observations above imply that there is a function
\begin{align*}
f:B_0&\longrightarrow B_{m-1}
\\
b&\longmapsto
\begin{cases}
b+m-1&\text{if }b<ym
\\
b-1&\text{if }b>ym
\end{cases}
\end{align*}
which is injective, and hence bijective. Now we partition each $B_g$ as $C_g\cup D_g$, where
\[
C_g=\lset{b\in B_g}{b<ym},\qquad D_g=\lset{b\in B_g}{b>ym}.
\]
The fact that $f$ is bijective together with observation \ref{underx} above means that
\[
C_g=C_0+g\text{ for }0\ls g\ls m-1,\tag{\tac}\label{lsy}
\]
so that in particular $\card{C_0}=\dots=\card{C_{m-1}}=r$, say. The fact that $f$ is bijective also gives
\[
D_{m-1}+1=D_0.\tag{\tac\tac}\label{gsy}
\]
So the combined bead moves from $D_{m-1},\dots,D_0$ have no effect on runner $0$. So (\ref{lsy}) implies that runner $0$ of $\ab\nu$ is obtained from runner $0$ of $\ab{\jo\la\pi}$ by adding a bead at position $b+m$ for all $b\in C_0$, and then removing a bead from position $b$ for all $b\in C_0$. In particular, there are exactly $a$ beads on runner $0$ in $\ab\nu$, and the last bead on runner $0$ is in position $ym$ or earlier. (\ref{lsy}) also means that all the positions before position $ym$ on runners other than runner $0$ are occupied in $\ab\nu$. This means that $\nu$ has the form $\jo\mu\rho$ with $\mu\in\pml$ and $\rho\in\calp$. Now we compare the beta-sets $\ber{a+1}\pi$ and $\ber{a+1}\rho$. If we let $A=\lset{c/m+1}{c\in C_0}$, then the above description of runner $0$ of $\ab\nu$ means that $\card A=r$, $A\cap\ber{a+1}\pi=\emptyset$, and $\ber{a+1}\rho=\ber{a+1}\pi\cup A\setminus(A-1)$. So $\pi\car\rho$ by \cref{carbeta}(2).

Another consequence of (\ref{lsy}) is that the combined bead moves from $C_{m-1},\dots,C_0$ have no effect on any runner other than runner $0$. This means that if we take $\ab\la$ and make the bead moves from $D_{m-1},\dots,D_0$ in turn, then we obtain $\ab\mu$. So $\mu$ appears in $F_{k-r}\la$ with non-zero coefficient.

Now we show that $\jo\mu\rho$ is $u$-\sep. Let $D=D_0\cup\dots\cup D_{m-1}$. Because $\pi\car\rho$, the last bead on runner $0$ in $\ab{\jo\mu\rho}$ is in position $(x+r)m$ or earlier. Let $d$ be the first empty position in $\ab{\jo\mu\rho}$ not on runner $0$. If $d>(x+k+u)m$, then certainly $\jo\mu\rho$ is $u$-\sep, so assume $d<(x+k+u)m$. Then position $d$ is occupied in $\ab{\jo\la\pi}$, so $d\in D$. Now observation \ref{underx} and (\ref{gsy}) imply that $d+1,d+2,\dots,(x+k+u)m$ all belong to $D$ as well. Hence
\[
m(k-r)=\card D>(x+k+u)m-d,
\]
giving $d>(x+r+u)m$, so that $\jo\mu\rho$ is $u$-\sep.

\smallskip
Now we need to prove the ``if'' part of the \lcnamecref{addrun1simple}. So suppose we are given $\mu,\rho,r$ satisfying the conditions in the \lcnamecref{addrun1simple}. Because $\mu$ appears in $F_{k-r}\la$ with non-zero coefficient, $\ab\mu$ is obtained from the \abd for $\la$ by making bead moves from $D_{m-1},\dots,D_0$ in turn, where $D_g\subset g+m\bbz$ and $\card{D_g}=k-r$ for each $g$. In addition, $D_0=D_{m-1}+1$ because $\la$ and $\mu$ have identical configurations on the $0$th runner. In the same way as in the last paragraph, we can show that every element of $D_0\cup\dots\cup D_{m-1}$ is greater than $(x+r+u)m$.

The assumption that $\pi\car\rho$ means (using \cref{carbeta}(2)) that $\ber{a+1}\rho=\ber{a+1}\pi\cup A\setminus(A-1)$ for some set $A$ with $\card A=r$ and $A\cap\ber{a+1}\pi=\emptyset$. So if we define
\[
C_g=\lset{(a-1)m+g}{a\in A}
\]
for each $g$, then $\ab{\jo\mu\rho}$ is obtained from $\ab{\jo\la\pi}$ by making bead moves from $C_{m-1}\cup D_{m-1},\dots,C_0\cup D_0$ in turn, so $\jo\mu\rho$ appears in $F_k(\jo\la\pi)$.
\end{pf}

Now we want to show that the coefficients appearing in \cref{addrun1simple} coincide.

\begin{propn}\label{addrun1}
Suppose $\la,\mu\in\pml$ and $\pi,\rho\in\calp$ with $\pi\car\rho$ for some $r$, and that $k\in\bbn$ is such that $k\gs r$ and $\jo\la\pi$ is $k$-\sep. Then
\[
\ip{F_k(\jo\la\pi)}{\jo\mu\rho}=\ip{F_{k-r}\la}{\mu}
\]
\end{propn}

\begin{pf}
\cref{addrun1simple} shows that one side of the given equation is non-zero \iff the other is, so we assume they are both non-zero. We use an abacus with $l+am+1$ beads, and we define the integer $y$ and the sets $B_g=C_g\cup D_g$ for $0\ls g\ls m-1$ as in the proof of \cref{addrun1simple}. Then $\mu$ is obtained from $\la$ by making bead moves from $D_{m-1},\dots,D_0$ in turn, and $\jo\mu\rho$ is obtained from $\jo\la\pi$ by making bead moves from $B_{m-1},\dots,B_0$ in turn. We remark that these sets are uniquely defined: if $\jo\mu\rho$ appears in $F_k(\jo\la\pi)=\f{-l}^{(k)}\f{1-l}^{(k)}\dots \f{-1-l}^{(k)}(\jo\la\pi)$, then we obtain $\jo\mu\rho$ by adding all the $(-1-l)$-nodes of $(\jo\mu\rho)\setminus(\jo\la\pi)$, then all the $(-2-l)$-nodes of $(\jo\mu\rho)\setminus(\jo\la\pi)$, and so on, with no choice at any stage.

So if we define partitions $\lb\la0,\dots,\lb\la m$ by setting $\lb\la0=\la$ and defining $\lb\la g$ from $\lb\la{g-1}$ by making bead moves from $D_{m-g}$ for $g\gs1$, then $\lb\la m=\mu$ and
\[
\ip{F_{k-r}\la}{\mu}=\prod_{g=1}^m\ip{\f{-g-l}^{(k-r)}\lb\la{g-1}}{\lb\la g},
\]
and the formula for the action of $\f{-g-l}^{(k-r)}$ says that
\[
\ip{\f{-g-l}^{(k-r)}\lb\la{g-1}}{\lb\la g}=q^{2n(\lb\la{g-1},\lb\la g)}.
\]
Similarly we define partitions $\lb\xi0,\dots,\lb\xi m$ by setting $\lb\xi0=\jo\la\pi$ and then making bead moves from $B_{m-1},\dots,B_0$, and derive a similar formula for $\ip{F_k(\jo\la\pi)}{\jo\mu\rho}$. We compare $n(\lb\la{g-1},\lb\la g)$ with $n(\lb\xi{g-1},\lb\xi g)$ for each $g$.

Taking $2\ls g\ls m-1$ first and letting $i=-g-l+m\bbz$, we compare the addable and removable $i$-nodes of $\lb\la{g-1},\lb\la g,\lb\xi{g-1},\lb\xi g$, by looking at runners $m-g$ and $m-g+1$ of their \abds. For any given row we may give the configuration of beads on the two positions on runner $m-g$ and $m-g+1$ in that row.

The rows on these runners where these four partitions differ are the rows containing elements of $C_{m-g}$, where
\begin{itemize}
\item
$\lb\la{g-1}$ and $\lb\la g$ both have $\miabb$,
\item
$\lb\xi{g-1}$ has $\miabn$,
\item
$\lb\xi g$ has $\mianb$,
\end{itemize}
and the rows containing elements of $D_{m-g}$, where
\begin{itemize}
\item
$\lb\la{g-1}$ and $\lb\xi{g-1}$ both have $\miabn$,
\item
$\lb\la g$ and $\lb\xi g$ both have $\mianb$.
\end{itemize}

We conclude that $\lb\la{g-1}$ and $\lb\xi{g-1}$ have exactly the same removable $i$-nodes while $\lb\la g$ and $\lb\xi g$ have exactly the same addable $i$-nodes. Moreover, none of the addable $i$-nodes of $\lb\xi g$ or removable $i$-nodes of $\lb\xi{g-1}$ lie to the left of any of the added nodes corresponding to the positions in $C_{m-g}$. So $n(\lb\la{g-1},\lb\la g)=n(\lb\xi{g-1},\lb\xi g)$, which means that $\ip{\f i^{(k-r)}\lb\la{g-1}}{\lb\la g}=\ip{\f i^{(k)}\lb\xi{g-1}}{\lb\xi g}$.

Now we consider the case $g=1$, by comparing runners $m-1$ and $0$ of $\lb\la0,\lb\la1,\lb\xi0,\lb\xi1$. To help in visualising addable and removable nodes, we imagine runner $m-1$ moved to lie to the left of runner $0$, and shifted down so that position $cm-1$ is directly to the left of position $cm$ for each $c$. So for a given partition when we show the abacus configuration on these runners in row $c$, we show positions $cm-1$ and $cm$. (We temporarily introduce a new position $-1$, which is taken to be occupied in all \abds.)

After position $ym$, the partitions $\lb\la0$ and $\lb\xi0$ agree, as do the partitions $\lb\la1$ and $\lb\xi1$. Now we look at the rows before position $ym$.

The \abds for $\lb\la0$ and $\lb\la1$ both have $\miabb$ in the first $a+1$ rows, and $\miabn$ in the remaining $y-a$ rows up to position $ym$.

Now consider the positions $cm-1$ and $cm$ in the \abds for $\lb\xi0$ and $\lb\xi1$, where $cm\ls y$. The configurations in these positions are determined by the sets $C_{m-1}$ and $\ber{a+1}\pi$:
\begin{itemize}
\item
if $cm-1\in C_{m-1}$, then $\lb\xi0$ has $\miabn$ and $\lb\xi1$ has $\mianb$;
\item
if $cm-1\notin C_{m-1}$ and $c\notin\ber{a+1}\pi$, then $\lb\xi0$ and $\lb\xi1$ both have $\miabn$;
\item
if $c\in\ber{a+1}\pi$, then $\lb\xi0$ and $\lb\xi1$ both have $\miabb$.
\end{itemize}
Note in particular that the second possibility happens exactly $y-a-r$ times. This enables us to compare $n(\lb\la0,\lb\la1)$ with $n(\lb\xi0,\lb\xi1)$ by examining addable and removable nodes. We find that
\[
n(\lb\xi0,\lb\xi1)=n(\lb\la0,\lb\la1)-(k-r)r+|A|,
\]
where
\[
A=\lset{(c,d)}{0\ls c<d,\ cm-1\notin C_{m-1},\ c\notin\ber{a+1}\pi,\ dm-1\in C_{m-1}}.
\]
We perform a similar analysis for the case $g=m$, and obtain
\[
n(\lb\xi{m-1},\lb\xi m)=n(\lb\la{m-1},\lb\la m)+(k-r)r-|B|,
\]
where
\[
B=\lset{0\ls c<d\ls y}{cm-1\notin C_{m-1},\ c\notin\ber{a+1}\pi,\ dm\in C_0}.
\]
Now recall from the proof of \cref{addrun1simple} that $C_{m-1}=C_0+m-1$. This then means that there is a bijection $A\to B$ given by $(c,d)\mapsto(c,d-1)$. We deduce that
\[
n(\lb\xi0,\lb\xi1)+n(\lb\xi{m-1},\lb\xi m)=n(\lb\la0,\lb\la1)+n(\lb\la{m-1},\lb\la m),
\]
and therefore
\[
\ip{\f{-l}^{(k)}\lb\xi0}{\lb\xi1}\ip{\f{1-l}^{(k)}\lb\xi{m-1}}{\lb\xi m}=\ip{\f{-l}^{(k-r)}\lb\la0}{\lb\la1}\ip{\f{1-l}^{(k-r)}\lb\la{m-1}}{\lb\la m},
\]
which is enough to complete the proof.
\end{pf}

\subsection{Bar-separated $h$-strict partitions}\label{barsepsec}

Now we prove similar results for $h$-strict partitions. Keep $l$ fixed, and suppose $\alpha$ is a \nice \bp of length $l$. Given a partition $\pi$, define $\jh\al\pi=\al\sqcup h\pi$. So the \abd for $\jh\al\pi$ with $h$ runners is obtained by taking the \abd for $\al$ and adding a bead at position $\pi_rh$ for each $r$. Let $b$ be the position of the last bead on runner $0$ (so $b=\pi_1h$), and let $f$ be the first empty position on any of runners $1,\dots,n$. Given $k\gs0$, say that $\jh\al\pi$ is \emph{$k$-\bsep} $f-b>kh$. (We write simply ``\bsep'' to mean ``$0$-\bsep''.)

Now for each $k\in\bbn$, define $\spin F_k=\spf0^{(2k)}\spf 1^{(2k)}\dots\spf{n-1}^{(2k)}\spf n^{(k)}$.
Then we can give an analogue of \cref{addrun1simple}.

\begin{propn}\label{addrun2simple}
Suppose $\al\in\pnice l$, $\pi\in\calp$ and $k,u\in\bbn_0$, and that $\jh\al\pi$ is $(k+u)$-\bsep. If $\ga\in\hstr$, then $\ip{\spin F_k(\jh\al\pi)}{\ga}\neq0$ \iff $\ga$ has the form $\jh\be\rho$, where
\begin{itemize}
\item
$\be\in\pnice l$ and $\rho\in\calp$,
\item
$\pi\car\rho$ for some $0\ls r\ls k$, and
\item
$\ip{\spin F_{k-r}\al}{\be}\neq0$.
\end{itemize}
Furthermore, if these conditions hold, then $\jh\be\rho$ is $u$-\bsep.
\end{propn}

\begin{pf}
We follow the structure of the proof of \cref{addrun1simple}, though the details here are slightly different. For any $h$-strict partition $\ga$, we write $\bab\ga$ for the \babd of $\ga$. As in the proof of \cref{addrun1simple}, given a finite set $B\subset\bbn_0$, we use the phrase \emph{making bead moves from $B$} to mean moving a bead from position $b$ to position $b+1$ for each $b\in B$.

Suppose $\ga$ appears in $\spin F_k(\jh\al\pi)$. Let $\be$ be the partition obtained from $\ga$ by removing all the positive parts divisible by $h$, and let $\rho$ be the partition obtained by taking all the parts of $\ga$ divisible by $h$ and dividing them all by $h$. The assumption that $\ip{\spin F_k(\jh\al\pi)}{\ga}\neq0$ means that we can obtain $\bab\ga$ from $\bab{\jh\al\pi}$ by making bead moves from
\[
B_n,B_{n+1}\cup B_{n-1},B_{n+2}\cup B_{n-2},\dots,B_{h-1}\cup B_0
\]
in turn, where each $B_g$ is a set of non-negative integers congruent to $g$ modulo $h$, and $\card{B_n}=k$ while $\card{B_{n+g}}+\card{B_{n-g}}=2k$ for $1\ls g\ls n$.

Let $x=\pi_1$. Then the lowest bead on runner $0$ in $\bab{\jh\al\pi}$ is in position $xh$. Because $\card{B_n}=k$, we can choose $y\in\{x,\dots,x+k\}$ such that $yh+n\notin B_n$. This has two consequences.
\begin{itemize}
\item
Because $\jh\al\pi$ is $k$-\bsep, positions $yh+1,\dots,yh+n$ are occupied in $\bab{\jh\al\pi}$. The assumption that $yh+n\notin B_n$ then means that $yh+g\notin B_g$ for $g=1,\dots,n-1$.
\item
Because $\al$ is \nice, runners $n+1,\dots,h-1$ are empty in $\bab{\jh\al\pi}$ and position $(y+1)h$ is unoccupied. The assumption that $yh+n\notin B_n$ then means that $yh+g\notin B_g$ for $g=n+1,\dots,h-1$, and $(y+1)h\notin B_0$.
\end{itemize}
Now write $B_g=C_g\cup D_g$ for each $g$, where
\[
C_g=\lset{b\in B_g}{b\ls yh},\qquad D_g=\lset{b\in B_g}{b>yh}.
\]
The fact that runners $n+1,\dots,h-1$ are empty in $\bab{\jh\al\pi}$, and that there are no beads on runner $0$ after position $yh$, means that $C_{n+g}\subseteq C_{n+g-1}+1$ and $D_{n+g}\subseteq D_{n+g-1}+1$ for $g=1,\dots,n$, and also that $D_0\subseteq D_n+n+1$. In addition, the fact that all the positions on runners $1,\dots,n$ before position $yh$ are occupied in $\bab{\jh\al\pi}$ means that $C_{n-g}\subseteq C_{n-g+1}-1$ for $g=1,\dots,n$. So
\[
\card{C_0}\ls\dots\ls\card{C_n}\gs\dots\gs\card{C_{h-1}}
\]
and
\[
\card{D_n}\gs\dots\gs\card{D_{h-1}}\gs\card{D_0}
\]

Now the fact that $\card{C_{h-1}}+\card{C_0}+\card{D_{h-1}}+\card{D_0}=2k=2\card{C_n}+2\card{D_n}$ gives equality everywhere. So $C_g=C_0+g$ for $g=1,\dots,h-1$ and $D_{n+g}=D_n+g$ for $g=1,\dots,n$, and $D_0=D_n+n+1$.

So the combined bead moves from $C_n,C_{n+1}\cup C_{n-1},\dots,C_{h-1}\cup C_0$ only affect runner $0$ in the \babd, where the effect is to add a bead at position $b+h$ for each $b\in C_0$ and then remove a bead from position $b$ for each $b\in C_0$. By \cref{carbeta}(1), this means that $\pi\car\rho$.

Now look at the other runners in $\bab\ga$. From what we have learned about the sets $C_g$ and $D_g$, we know that there are $l$ beads on runners $1,\dots,n$, and no beads on runners $n+1,\dots,h-1$. So $\be\in\pnice l$, and $\ga=\jh\be\rho$. Furthermore, $\bab\be$ is obtained from $\bab\al$ by making bead moves from $D_n,D_{n+1}\cup D_{n-1},\dots,D_{h-1}\cup D_0$, and $\card{D_g}=k-r$ for each $g$, so $\ip{\spin F_{k-r}\al}{\be}\neq0$.

We also need to show that $\jh\be\rho$ is $u$-\bsep; this is done in the same way as the corresponding part of \cref{addrun1simple}.

\smallskip
Now we prove the ``if'' part of the \lcnamecref{addrun2simple}. Suppose we are given $\be,\rho,r$ as in the \lcnamecref{addrun2simple}. Because $\ip{\spin F_{k-r}\al}{\be}\neq0$, there are sets $D_0,\dots,D_{h-1}\subseteq\bbn_0$ with $D_g\subset g+h\bbz$ for each $g$ and $\card{D_{n+g}}+\card{D_{n-g}}=2(k-r)$ for $g=0,\dots,n$, such that $\bab\be$ is obtained from $\bab\al$ by making bead moves from $D_n,D_{n+1}\cup D_{n-1},\dots,D_{h-1}\cup D_0$ in turn. Because $\al$ and $\be$ are both \nice, we obtain $D_{n+g}=D_n+g$ for $g=1,\dots,n$, so that in fact $\card{D_g}=k-r$ for every $g$.

The assumption that $\pi\car\rho$ means that there is a set $A\subset\bbn$ with $\card A=r$ such that $\rho$ is obtained from $\pi$ by adding a part equal to $a$ for every $a\in A$, and then removing a part equal to $a-1$ for each $a\in A$. So if we define
\[
C_g=\lset{(a-1)h+g}{a\in A},\qquad B_g=C_g\cup D_g
\]
for each $g$, then $\bab{\jh\be\rho}$ is obtained from $\bab{\jh\al\pi}$ by making bead moves from $B_n,B_{n+1}\cup B_{n-1},B_{h-1}\cup B_0$ in turn, so $\jh\be\rho$ appears in $\spin F_k(\jh\al\pi)$.
\end{pf}

Now we compare the coefficients appearing in \cref{addrun2simple}. This is similar to \cref{addrun1}, though here the statement is more complicated. Given partitions $\pi,\rho$ with $\pi\car\rho$, we let $\hsb\pi\rho$ be defined as in \cref{newpierisec}, \emph{with the indeterminate $t$ replaced by $-q^2$}; that is,
\[
\hsb\pi\rho=\prod_{\substack{c\gs1\\\rho'_c=\pi'_c\\\rho'_{c+1}>\pi'_{c+1}}}(1-(-q^2)^{\pi'_c-\pi'_{c+1}}).
\]

\begin{propn}\label{addrun2}
Suppose $\al,\be\in\pnice l$ and $\pi,\rho\in\calp$ with $\pi\car\rho$ for some $r$, and that $\jh\al\pi$ is $k$-\bsep for some $k\gs r$. Then
\[
\ip{\spin F_k(\jh\al\pi)}{\jh\be\rho}=\hsb\pi\rho\ip{\spin F_{k-r}\al}{\be}.
\]
\end{propn}

\begin{pf}
\cref{addrun2simple} shows that one side of the given equation is non-zero \iff the other is, so we assume they are both non-zero. We define the integer $y$ and the sets $B_g=C_g\cup D_g$ for $0\ls g\ls h-1$ as in the proof of \cref{addrun2simple}. Then $\be$ is obtained from $\al$ by making bead moves from $D_n,D_{n+1}\cup D_{n-1},\dots,D_{h-1}\cup D_0$ in turn, and $\jh\be\rho$ is obtained from $\jh\al\pi$ by making bead moves from $B_n,B_{n+1}\cup B_{n-1},\dots,B_{h-1}\cup B_0$ in turn. As in the proof of \cref{addrun1}, these sets are uniquely defined.

Define partitions $\lb\al0,\dots,\lb\al m$ by letting $\lb\al0=\al$, and then constructing $\lb\al{g+1}$ from $\lb\al g$ by making bead moves from $B_{n+g}\cup B_{n-g}$, for $g=0,\dots,n-1$. Then $\lb\al m=\be$, and
\begin{alignat*}2
\ip{\spin F_{k-r}\al}{\be}&=\ip{\spf n^{(k-r)}\al}{\lb\al1}\ \mathrlap{\times\ \prod_{g=1}^n\ip{\spf{n-g}^{(2(k-r))}\lb\al g}{\lb\al{g+1}}.}
\\
\intertext{Correspondingly,}
\ip{F_k(\jh\al\pi)}{\jh\be\rho}&=\ip{\spf n^{(k)}(\jh\al\pi)}{\jh{\lb\al1}\pi}\ &&\times\ \prod_{g=1}^{n-1}\ip{\spf{n-g}^{(2k)}(\jh{\lb\al g}\pi)}{\jh{\lb\al{g+1}}\pi}
\\
&&&\times\ \ip{\spf0^{(2k)}(\jh{\lb\al n}\pi)}{\jh\be\rho}.
\end{alignat*}
First we compare the terms
\[
\ip{\spf n^{(k-r)}\al}{\lb\al1}=q^{4\spin n(\al,\lb\al1)},
\qquad
\ip{\spf n^{(k)}(\jh\al\pi)}{\jh{\lb\al1}\pi}=q^{4\spin n(\jh\al\pi,\jh{\lb\al1}\pi)}.
\]
We can compute the integers $\spin n(\al,\lb\al1)$ and $\spin n(\jh\al\pi,\jh{\lb\al1}\pi)$ by looking at runners $n$ and $n+1$ in the bar-abacus. As in the proof of \cref{addrun1}, we may draw a small diagram showing the configuration of beads on these runners in a given row. These two runners are identical in the four partitions $\al,\lb\al1,\jh\al\pi,\jh{\lb\al1}\pi$, except:
\begin{itemize}
\item
in the positions corresponding to elements of $C_n$, where $\al$, $\lb\al1$ and $\jh\al\pi$ all have $\mibbn$ while $\jh{\lb\al1}\pi$ has $\mibnb$;
\item
in the positions corresponding to elements of $D_n$, where $\al$ and $\jh\al\pi$ have $\mibbn$ while $\lb\al1$ and $\jh{\lb\al1}\pi$ have $\mibnb$.
\end{itemize}

Now comparing the calculation of $\spin n(\al,\lb\al1)$ and $\spin n(\jh\al\pi,\jh{\lb\al1}\pi)$, we obtain
\[
\spin n(\jh\al\pi,\jh{\lb\al1}\pi)=\spin n(\al,\lb\al1)-(k-r)r+\card A,
\]
where
\[
A=\lset{0\ls c<d}{ch+n\notin C_n\ni dh+n}.
\]

Next we take $1\ls g\ls n-1$, and compare the terms
\[
\ip{\spf{n-g}^{(2(k-r))}\lb\al g}{\lb\al{g+1}}=q^{2\spin n(\lb\al g,\lb\al{g+1})},\qquad\ip{\spf{n-g}^{(2k)}(\jh{\lb\al g}\pi)}{\jh{\lb\al{g+1}}\pi}=q^{2\spin n(\jh{\lb\al g}\pi,\jh{\lb\al{g+1}}\pi)}.
\]
We can compare the integers $\spin n(\lb\al g,\lb\al{g+1})$ and $\spin n(\jh{\lb\al g}\pi,\jh{\lb\al{g+1}}\pi)$ by looking at runners $n-g,n-g+1,n+g,n+g+1$ of the bar-abacus.

On runners $n-g$ and $n-g+1$ the four partitions have identical configurations, except:
\begin{itemize}
\item
in the positions corresponding to elements of $C_{n-g}$, where $\lb\al g$ and $\lb\al{g+1}$ have $\mibbb$, while $\jh{\lb\al g}\pi$ has $\mibbn$ and $\jh{\lb\al{g+1}}\pi$ has $\mibnb$;
\item
in the positions corresponding to elements of $D_{n-g}$, where $\lb\al g$ and $\jh{\lb\al g}\pi$ have $\mibbn$ while $\lb\al{g+1}$ and $\jh{\lb\al{g+1}}\pi$ have $\mibnb$.
\end{itemize}
A similar statement holds for runners $n+g$ and $n+g+1$, but with $\mibnn$ in place of $\mibbb$. Now comparing the calculation of $\spin n(\lb\al g,\lb\al{g+1})$ and $\spin n(\jh{\lb\al g}\pi,\jh{\lb\al{g+1}}\pi)$ using the bar-abacus gives
\[
\spin n(\lb\al g,\lb\al{g+1})=\spin n(\jh{\lb\al g}\pi,\jh{\lb\al{g+1}}\pi).
\]

Finally we look at the terms
\[
\ip{\spf0^{(2(k-r))}\lb\al n}{\be}=N_1q^{\spin n(\lb\al n,\be)},\qquad\ip{\spf0^{(2k)}(\jh{\lb\al n}\pi)}{\jh\be\rho}=N_2q^{\spin n(\jh{\lb\al n}\pi,\jh\be\rho)}.
\]
Here $N_1$ and $N_2$ are products of terms $1-(-q^2)^b$ which we will deal with at the end of the proof.

To compare the calculation of $\spin n(\lb\al n,\be)$ and $\spin n(\jh{\lb\al n}\pi,\jh\be\rho)$ it is easier to use Young diagrams directly rather than the bar-abacus. Given a pair of columns $c<d$ where $c,d$ are each congruent to $0$ or $1$ modulo $h$, let's write
\[
\ar_{c,d}(\lb\al n,\be)=
\begin{cases}
1&\parbox{205pt}{if there is a node of $\be\setminus\lb\al n$ in column $d$ and a $0$-bar-addable node of $\be$ in column $c$}
\\[9pt]
-1&\parbox{205pt}{if there is a node of $\be\setminus\lb\al n$ in column $d$ and a $0$-bar-removable node of $\lb\al n$ in column $c$}
\\[9pt]
0&\text{otherwise}.
\end{cases}
\]
We define $\ar_{c,d}(\jh{\lb\al n}\pi,\jh\be\rho)$ similarly. Then
\[
\spin n(\jh{\lb\al n}\pi,\jh\be\rho)-\spin n(\lb\al n,\be)=\sum_{c<d}\left(\ar_{c,d}(\jh{\lb\al n}\pi,\jh\be\rho)-\ar_{c,d}(\lb\al n,\be)\right).
\]
We can explicitly list the pairs $c<d$ for which $\ar_{c,d}(\jh{\lb\al n}\pi,\jh\be\rho)\neq\ar_{c,d}(\lb\al n,\be)$. These occur in two ways.

\begin{itemize}
\item
If $d-1\in C_{h-1}\cup C_0$ while $c-1\notin C_{h-1}\cup C_0$, then $\ar_{c,d}(\lb\al n,\be)=0$ while $\ar_{c,d}(\jh{\lb\al n}\pi,\jh\be\rho)=-1$.
\item
If $d-1\in D_{h-1}\cup D_0$ and $c-1\in C_{h-1}\cup C_0$, then $\ar_{c,d}(\lb\al n,\be)=-1$ while $\ar_{c,d}(\jh{\lb\al n}\pi,\jh\be\rho)=0$.
\end{itemize}
Summing over all these pairs, we obtain
\begin{align*}
\spin n(\jh{\lb\al n}\pi,\jh\be\rho)-\spin n(\lb\al n,\be)&=(\card{C_0}+\card{C_{h-1}})(\card{D_0}+\card{D_{h-1}})-\card B
\\
&=4r(k-r)-\card B,
\end{align*}
where
\[
B=\lset{(c,d)}{1\ls c<d,\ c\equiv0,1\ppmod h,\ c-1\notin(C_0\cup C_{h-1}),\ d-1\in(C_0\cup C_{h-1})}.
\]
But now recalling that $C_0=C_n-n$ and $C_{h-1}=C_n+n$ gives
\[
B=\lset{(ch\pm n,dh\pm n)}{(c,d)\in A},
\]
where $A$ is the set defined earlier in the proof. So $\card B=4\card A$.

Now combining all the cases we have computed to calculate $\ip{F_k(\jh\al\pi)}{\jh\be\rho}/\ip{\spin F_{k-r}\al}{\be}$, the powers of $q$ cancel to give
\[
\frac{\ip{F_k(\jh\al\pi)}{\jh\be\rho}}{\ip{\spin F_{k-r}\al}{\be}}=\frac{N_2}{N_1},
\]
and we just need to show that $N_2/N_1=\hsb\pi\rho$. Recall that $N_1$ is the product of terms $1-(-q^2)^b$, where we take the product over all $c\gs1$ such that $\be\setminus\lb\al n$ includes a node in column $ch+1$ but not column $c$, and $b$ is the number of times $c$ appears in $\lb\al n$. But in fact there are no such $c$, because the nodes of $\be\setminus\lb\al n$ correspond to elements of $D_{h-1}\cup D_0$, and $D_0=D_{h-1}+1$. So $N_1=1$. So we need to consider $N_2$, which is defined in a similar way. The values of $c$ for which $\jh\be\rho\setminus\jh{\lb\al n}\pi$ includes a node in column $ch+1$ but not column $c$ are the values such that $ch\in C_0$ but $ch-1\notin C_{h-1}$. These in turn are precisely the values of $c$ for which $\rho'_{c+1}>\pi'_{c+1}$ while $\rho'_c=\pi'_c$. As a result, $N_2=\hsb\pi\rho$, as required.
\end{pf}

\section{Comparing canonical bases}\label{comparingsec}

In this section we give our main results comparing canonical bases. In order to use the results of \cref{sepsec,barsepsec}, we need to compare the actions of the operators $F_k$ and $\spin F_k$ defined in \cref{nextbitsec}. We begin with the following \lcnamecref{samecoeff}.

\begin{propn}\label{samecoeff}
Suppose $\al,\be\in\pnice l$. Then
\[
\ip{\spin F_k\al}{\be}=\ip{F_k\psi\al}{\psi\be}.
\]
\end{propn}

\begin{pf}
For this proof we write $\la=\psi\al$ and $\mu=\psi\be$. First we show that one side of the equation is non-zero \iff the other is. We use the terminology and notation relating to the (bar)-abacus from the proofs of \cref{addrun1simple,addrun1,addrun2simple,addrun2}; we use \abds with $l+1$ beads for $\la$ and $\mu$. We define functions
\begin{alignat*}2
\hf:\bbz&\longrightarrow\bbz&\qquad\fh:\bbz&\longrightarrow\bbz
\\
b&\longmapsto b+(m-h)\inp{\frac bh}&b&\longmapsto b+(h-m)\inp{\frac bm}.
\end{alignat*}
Observe that $\hf$ and $\fh$ restrict to mutually inverse bijections between $g+h\bbz$ and $g+m\bbz$ for each $0\ls g\ls n$, and that position $b$ is occupied in $\ab\la$ \iff position $\fh(b)$ is occupied in $\bab\al$ (and similarly for $\mu$ and $\be$).

Suppose first that $\ip{\spin F_k\al}\be\neq0$. Then there are sets $D_0,\dots,D_{h-1}\subset\bbn_0$ with $D_g\subset g+h\bbz$ and $\card{D_{n-g}}+\card{D_{n+g}}=2k$ for each $g$, such that $\bab\be$ is obtained from $\bab\al$ by making bead moves from $D_n,D_{n+1}\cup D_{n-1},\dots,D_{h+1}\cup D_0$ in turn. The assumption that $\al,\be\in\pnice l$ implies that $D_{n+g}=D_{n+g-1}+1$ for $g=1,\dots,n$, and that $D_0=D_{h-1}+1$. So in fact $\card{D_g}=k$ for each $g$. Now if we define $E_g=\hf(D_g)$ for $g=0,\dots,n$, then $E_g\subset g+m\bbz$ and $\card{E_g}=k$ for each $g$, and $\ab\mu$ can be obtained from $\ab\la$ by making bead moves from $E_n,E_{n-1},\dots,E_0$ in turn. So $\ip{F_k\la}\mu\neq0$.

The converse is very similar: if $\ip{F_k\la}\mu\neq0$, then there are sets $E_0,\dots,E_n$ with the properties given above and satisfying $E_0=E_n+1$. Now if we define $D_g=\fh(E_g)$ for $g=0,\dots,n$, and $D_g=D_{g-1}+1$ for $g=n+1,\dots,h-1$, then the sets $D_g$ have the properties in the above paragraph and $\bab\be$ is obtained from $\bab\al$ by making bead moves from $D_n,D_{n+1}\cup D_{n-1},\dots,D_{h-1}\cup D_0$. So $\ip{\spin F_k\al}\be\neq0$.

So we assume that both sides of the equation are non-zero, and define the sets $D_0,\dots,D_{h-1}$ and $E_0,\dots,E_n$ as in the last two paragraphs. Also, for $g=1,\dots,n$ let $F_g$ be the set of parts $b$ of $\al$ congruent to $g$ modulo $h$ such that $b\notin D_g$; in other words, $F_g$ is the set of positions of beads on runner $g$ of $\bab\al$ that are \emph{not} moved in constructing $\be$ from $\al$. We define the partitions $\lb\al0,\dots,\lb\al m$ as in the proof of \cref{addrun2}, and the partitions $\lb\la0,\dots,\lb\la m$ as in the proof of \cref{addrun1} (using $E_0,\dots,E_n$ in place of $D_0,\dots,D_n$) , so that
\begin{align*}
\ip{\spin F_k\al}{\be}&=\ip{\spf n^{(k)}\al}{\lb\al1}\ \times\ \prod_{g=1}^n\ip{\spf{n-g}^{(2k)}\lb\al g}{\lb\al{g+1}},
\\
\ip{F_k\la}{\mu}&=\prod_{g=0}^n\ip{\f{-g-1-l}^{(k)}\lb\la g}{\lb\la{g+1}}.
\end{align*}
We will compare these two expressions term by term. First consider $\ip{\spf n^{(k)}\al}{\lb\al1}=q^{4\spin n(\al,\lb\al1)}$ and $\ip{\f{-1-l}^{(k)}\la}{\lb\la1}=q^{2n(\la,\lb\la1)}$. Since runner $n+1$ of $\bab\al$ is empty, $\al$ has no $n$-bar-removable nodes, so
\[
\spin n(\al,\lb\al1)=\card{\lset{(c,d)}{c<d,\ c\in F_n,\ d\in D_n}}.
\]
A similar statement applies for $n(\la,\lb\la1)$, and the relationship between the \abds then gives $n(\la,\lb\la1)=\spin n(\al,\lb\al1)$.

Now take $1\ls g\ls n-1$, and consider $\ip{\spf{n-g}^{(2k)}\lb\al g}{\lb\al{g+1}}=q^{2\spin n(\lb\al g,\lb\al{g+1})}$ and $\ip{\f{-g-1-l}^{(k)}\lb\la g}{\lb\la{g+1}}=q^{2n(\lb\la g,\lb\la{g+1})}$. By considering $\bab{\lb\al g}$, we obtain
\begin{align*}
\spin n(\lb\al g,\lb\al{g+1})=&\cardx{\rset{(c,d)}{c<d,\ c\in F_{n-g},\ d\in D_{n-g}\cup D_{n+g}}}
\\
&-\cardx{\rset{(c,d)}{c<d,\ c\in F_{n-g+1},\ d\in D_{n-g}\cup D_{n+g}}}.
\intertext{(Note that the terms in the minuend with $c\in F_{n-g}\cap(F_{n-g+1}-1)$ do not correspond to bar-addable nodes. But these terms cancel with the terms in the subtrahend with $c\in F_{n-g+1}\cap(F_{n-g}+1)$, which do not correspond to bar-removable nodes.) Doing the same calculation for $\lb\la g$ and applying $\fh$,}
n(\lb\la g,\lb\la{g+1})=&\cardx{\rset{(c,d)}{c<d,\ c\in F_{n-g},\ d\in D_{n-g}}}
\\
&-\cardx{\rset{(c,d)}{c<d,\ c\in F_{n-g+1},\ d\in D_{n-g}}}.
\end{align*}
Finally, we look at $\ip{\spf0^{(2k)}\lb\al n}{\be}=q^{\spin n(\lb\al n,\be)}$ and $\ip{\f{-l}^{(k)}\lb\la n}{\mu}=q^{2n(\lb\la n,\be)}$. Now there are two added nodes of bar-residue $0$ corresponding to each element of $D_{h-1}$, two $0$-bar-removable nodes corresponding to each element of $F_1\setminus\{1\}$, and a $0$-bar-addable node in column $1$ provided $1\notin F_1$. So
\begin{align*}
\spin n(\lb\al n,\be)&=2k-4\cardx{\lset{(c,d)}{c<d,\ c\in F_1,\ d\in D_{h-1}}}.
\\
\intertext{Correspondingly,}
n(\lb\la n,\mu)&=k-\cardx{\lset{(c,d)}{c<d,\ c\in F_1,\ d\in D_{h-1}}}.
\end{align*}
Combining these cases, we find that $\ip{\spin F_k\al}{\be}=q^{2s}\ip{F_k\la}{\mu}$, where
\[
s=\sum_{g=1}^n\cardx{\rset{(c,d)}{c<d,\ c\in F_{n-g+1},\ d\in D_{n+g-1}}}-\cardx{\rset{(c,d)}{c<d,\ c\in F_{n-g+1},\ d\in D_{n+g}}}.
\]
But $D_{n+g}=D_{n+g-1}+1$ for each $g$, so each summand is zero, and hence $\ip{\spin F_k\al}{\be}=\ip{F_k\la}{\mu}$.
\end{pf}

Now we are in a position to prove our main result comparing canonical basis elements, extending \cref{firstmain}. We consider the two sets of partitions
\[
\lset{\jh\al\pi}{\al\in\pnice l,\ \pi\in\calp},\qquad\lset{\jo\la\pi}{\la\in\pml,\ \pi\in\calp}.
\]
If we let $\psi$ be the bijection from \cref{firstcomparesec}, then there is a bijection between these two sets given by $\jh\al\pi\mapsto\jo{\psi\al}\pi$. Our aim is to extend the linear map $\Psi$ to give a correspondence between the canonical basis elements. The obvious map $\jh\la\pi\mapsto\jo{\psi\al}\pi$ does not work, because of the factor $\hsb\pi\rho$ appearing in \cref{addrun2}.

In fact, we need to restrict quite substantially the set of partitions that we consider. Suppose $\la\in\pml$ and $\pi\in\calp$. Given $u\gs0$, say that the partition $\jo\la\pi$ is \emph{$u$-\ssep} if $\jo\la\rho$ is $u$-\sep for every partition $\rho$ with $\card\rho=\card\pi$. (Equivalently, $\jo\la\pi$ is $u$-\ssep if it is $(u+\card\pi-\pi_1)$-\sep.)

Similarly, given $\al\in\pnice l$ and $\pi\in\calp$, say that $\jh\al\pi$ is \emph{$u$-\sbsep} if $\jh\al\rho$ is $u$-\bsep for every $\rho$ with $\card\rho=\card\pi$. (Equivalently, $\jh\al\pi$ is $u$-\sbsep \iff $\jo{\psi\al}\pi$ is $u$-\ssep.)

For the case $u=0$, we will simply say ``\ssep'' to mean ``$0$-\ssep'', and similarly for \sbsep.

Our next step is to show how the $u$-\ssep condition is affected when we apply the operators $F_k$. The following is a variation on the final statement of \cref{addrun1simple}.

\begin{lemma}\label{fssep}
Suppose $\la,\mu\in\pml$ and $\pi,\rho\in\calp$, and that $k>0$ with $\ip{F_k(\jo\la\pi)}{\jo\mu\rho}\neq0$. If $\jo\la\pi$ is $(k+u)$-\ssep for some $u\gs0$, then $\jo\mu\rho$ is $u$-\ssep.
\end{lemma}

\begin{pf}
The statement that $\jo\la\pi$ is $(k+u)$-\ssep is the same as saying that $\jo\la\pi$ is $(k+u+\card\pi-\pi_1)$-\sep. Since $\jo\mu\rho$ appears with non-zero coefficient in $F_k(\jo\la\pi)$, \cref{addrun1simple} shows that $\jo\mu\si$ does as well, where $\si$ is the partition obtained by adding $\card\rho-\card\pi$ nodes to $\pi$ at the end of the first row. So from the final statement of \cref{addrun1simple}, $\jo\mu\si$ is $(u+\card\pi-\pi_1)$-\sep, which is the same as $(u+\card\si-\si_1)$-\sep, so $\jo\mu\si$ is $u$-\ssep, and therefore $\jo\mu\rho$ is $u$-\ssep.
\end{pf}

Now for each $u$ define
\begin{align*}
\sspu u&=\lspan{\jo\la\pi}{\la\in\pml,\ \pi\in\calp,\ \jo\la\pi\text{ $u$-\ssep}},
\\
\bsspu u&=\lspan{\jh\al\pi}{\al\in\pnice l,\ \pi\in\calp,\ \jh\al\pi\text{ $u$-\sbsep}}.
\end{align*}
We write $\ssp$ and $\bssp$ for $\sspu0$ and $\bsspu0$.

\cref{fssep} shows that $F_k$ maps $\sspu{k+u}$ to $\sspu u$. Similarly, $\spin F_k$ maps $\bsspu{k+u}$ to $\bsspu u$. To extend \cref{samecoeff} to connect the actions of $F_k$ and $\spin F_k$, we need to define a linear bijection between $\bssp$ and $\ssp$. Recalling the Kostka polynomials $K_{\rho\pi}(t)$ from \cref{symfnsec}, we define
\[
\asa:\bssp\longrightarrow\ssp
\]
by setting
\[
\asa(\jh\al\pi)=\sum_{\rho\in\calp}K_{\rho\pi}(-q^2)\jo{\psi\al}\rho
\]
and extending linearly. ($\asa$ is bijective because the matrix of Kostka polynomials $K_{\rho\pi}(t)$ for $\rho,\pi\in\calp(r)$ is invertible for each $r$.)

Now we can reconcile \cref{addrun1,addrun2} via the following \lcnamecref{sasfk}.

\begin{propn}\label{sasfk}
Suppose $v\in\bsspu k$. Then $\asa(\spin F_kv)=F_k\asa(v)$.
\end{propn}

\begin{pf}
By linearity we can assume $v=\jh\al\pi$ for some $k$-\sbsep partition $\jh\al\pi$. Since both sides of the equation lie in $\ssp$, it suffices to show that
\[
\ip{\asa(\spin F_k(\jh\al\pi))}{\jo\mu\si}=\ip{F_k\asa(\jh\al\pi)}{\jo\mu\si}
\]
for each \sep partition $\jo\mu\si$.

We use the results of \cref{symfnsec}, with the variable $t$ specialised to $-q^2$. So we let $P_\la,Q_\la$ be the two types of Hall--Littlewood symmetric functions defined for $t=-q^2$, and $\lan\,,\,\ran$ the Hall--Littlewood inner product.

\allowdisplaybreaks
Now
\begin{align*}
\ip{\asa(\spin F_k(\jh\al\pi))}{\jo\mu\si}
&=\sum_{\substack{\tau\in\calp\\\pi\car\tau}}\sum_{\be\in\pnice l}\hsb\pi\tau\ip{\spin F_k\al}{\be}\ip{\asa\left(\jh\be\tau\right)}{\jo\mu\si}
\tag*{(by \cref{addrun2})}
\\
&=\sum_{\substack{\tau\in\calp\\\pi\car\tau}}K_{\si\tau}(-q^2)\hsb\pi\tau\ip{\spin F_k\al}{\psi^{-1}\mu}
\tag*{(from the definition of $\asa$)}
\\
&=\sum_{\substack{\tau\in\calp\\\pi\car\tau}}\lan s_\si,Q_\tau\ran\hsb\pi\tau\ip{F_k\psi\al}{\mu}
\tag*{(by \cref{samecoeff} and the definition of the coefficients $K_{\si\tau}(t)$)}
\\
&=\sum_{\tau\in\calp}\lan s_\si,Q_\tau\ran\lan\partial_{(r)}P_\tau,Q_\pi\ran\ip{F_k\psi\al}{\mu}
\tag*{(by \cref{dualpieri})}
\\
&=\lan\partial_{(r)}s_\si,Q_\pi\ran\ip{F_k\psi\al}{\mu}
\tag*{(since $s_\si=\sum_\tau\lan s_\si,Q_\tau\ran P_\tau$)}
\\
&=\sum_{\substack{\rho\in\calp\\\rho\car\si}}\lan s_\rho,Q_\pi\ran\ip{F_k\psi\al}{\mu}
\tag*{(by the Pieri rule for Schur functions)}
\\
&=\sum_{\rho\in\calp}K_{\rho\pi}(-q^2)\ip{F_k\jo{\psi\al}\rho}{\jo\mu\si}
\tag*{(by \cref{addrun1})}
\\
&=\ip{F_k\asa(\jh\al\pi)}{\jo\mu\si}
\tag*{(from the definition of $\asa$).\qedhere}
\end{align*}

\end{pf}

Now we want to use \cref{sasfk} to compare canonical basis elements. Given partitions $\pi,\rho$, let's write $\rho\succcurlyeq\pi$ if either $\card\rho<\card\pi$ or $\rho\dom\pi$. Part (2) of the next \lcnamecref{sscbv} is our main result comparing canonical bases.

\begin{thm}\label{sscbv}
Suppose $\al\in\pnice l$ and $\pi\in\calp$, and that $\jh\al\pi$ is $u$-\sbsep for some $u\gs0$. Then:
\begin{enumerate}
\item
$\scb h(\jh\al\pi)$ is a linear combination of basis elements $\jh\be\rho$ in which $\rho\succcurlyeq\pi$ and $\jh\be\rho$ is $u$-\sbsep;
\item
$\asa(\scb h(\jh\al\pi))=\cb m(\jo{\psi\al}\pi)$.
\end{enumerate}
\end{thm}

\newcommand\vc V

\begin{pf}
We proceed by induction on $\pi$, using the order $\succcurlyeq$. The case $\pi=\vn$ follows from \cref{firstmain}, so assume $\pi\neq\vn$, and that the theorem is true if $\pi$ is replaced by any $\rho$ with $\rho\succ\pi$.

Let $k$ be the last non-zero part of $\pi$, let $\pi^-$ denote the partition obtained by deleting this last part, and consider the bar-invariant vector $\vc=\spin F_k\scb h(\jh\al{\pi^-})$. Since $\jh\al{\pi^-}$ is $(k+u)$-\sbsep, the inductive hypothesis says that $\scb h(\jh\al{\pi^-})$ is a linear combination of basis elements $\jh\ga\si$ for which $\si\succ\pi^-$ and $\jh\ga\si$ is $(k+u)$-\sbsep. Now \cref{addrun2simple} shows that the terms $\jh\be\rho$ appearing in $\vc$ all have the properties that $\si\car\rho$ for some $r\ls k$ and $\jh\be\rho$ is $u$-\sbsep. If $\si\car\rho$ with either $\card\si<\card{\pi^-}$ or $r<k$, then $\card\rho<\card\pi$.  On the other hand, if $\si\dom\pi^-$ and $\si\cak\rho$, then $\rho\dom\pi$ by \cref{carlem}, with equality only if $\si=\pi^-$.

So $\vc$ is a linear combination of terms $\jh\be\rho$ with $\rho\succcurlyeq\pi$. Moreover, $\hsb{\pi^-}\pi=1$, so for any $\be$ the coefficient of $\jh\be\pi$ in $\vc$ is the same as the coefficient of $\jh\be{\pi^-}$ in $\scb h(\jh\al{\pi^-})$. This means in particular that $\jh\al\pi$ occurs with coefficient $1$, and that for any $\be\neq\al$ the coefficient of $\jh\be\pi$ in $\vc$ is divisible by $q$.

Because $\vc$ is bar-invariant, we can write it as a linear combination of canonical basis vectors. The properties of $\vc$ described in the last paragraph and the inductive hypothesis mean that when we do this, the canonical basis vector $\scb h(\jh\al\pi)$ must occur with coefficient $1$. So we can write
\[
V=\scb h(\jh\al\pi)+\sum_{\be,\rho}a_{\be\rho}\scb h(\jh\be\rho),\tag{\tac}
\]
with the sum being over pairs $\be,\rho$ for which $\rho\succ\pi$ and $\jh\be\rho$ is $k$-\sbsep, and each coefficient $a_{\be\rho}$ lies in $\bbc(q+q^{-1})$. Now property (1) for $\scb h(\jh\al\pi)$ follows from the inductive hypothesis and the corresponding property for $V$ shown above.

Applying $\asa$ to both sides of (\tac) and rearranging, we obtain
\begin{align*}
\asa\left(\scb h(\jh\al\pi)\right)
&=\asa\left(\spin F_k\scb h(\jh\al{\pi^-})\right)-\sum_{\be,\rho}a_{\be\rho}\asa\left(\scb h(\jh\be\rho)\right)
\\
&=F_k\asa\left(\scb h(\jh\al{\pi^-})\right)-\sum_{\be,\rho}a_{\be\rho}\asa\left(\scb h(\jh\be\rho)\right)
\tag{by \cref{sasfk}}\\
&=F_k\cb m(\jh{\psi\al}{\pi^-})-\sum_{\be,\rho}a_{\be\rho}\cb m(\jo{\psi\be}\rho)
\tag*{(by the inductive hypothesis).}
\end{align*}
This means in particular that $\asa\left(\scb h(\jh\al\pi)\right)$ is bar-invariant. Moreover, $\asa\left(\scb h(\jh\al\pi)\right)$ equals $\jh\al\pi$ plus a linear combination of other standard basis vectors with coefficients divisible by $q$; this comes from the corresponding property for $\scb h(\jh\al\pi)$, the definition of $\asa$ and the fact (\cref{kostkaprops}) that the inverse Kostka matrix $K^{-1}(-q^2)$ is unitriangular with off-diagonal entries divisible by $q^2$. So the defining properties of the canonical basis mean that $\asa\left(\scb h(\jh\al\pi)\right)=\cb m(\jo{\psi\al}\pi)$.
\end{pf}

\begin{eg}
We will be most interested in applying \cref{sscbv} in weight spaces which are completely contained within $\bssp$, so we give an example to illustrate this situation. Take $m=3$, and consider the weight space in $\fo$ spanned by partitions with $3$-core $(2^2,1^2)$ and $3$-weight $3$. The canonical basis coefficients are given in \cref{examp1}, where the $(\la,\mu)$-entry of the given matrix is~$\dn\la\mu$).
\begin{figure}[t]\footnotesize
\[
\begin{array}{c|cccccccccc|}
&\rotatebox{90}{$(2^2,1^{11})$}
&\rotatebox{90}{$(2^5,1^5)$}
&\rotatebox{90}{$(3^3,2^2,1^2)$}
&\rotatebox{90}{$(4,3,1^8)$}
&\rotatebox{90}{$(4,3,2^3,1^2)$}
&\rotatebox{90}{$(4,3^2,2,1^3)$}
&\rotatebox{90}{$(4,3^2,2^2,1)$}
&\rotatebox{90}{$(5,3^2,1^4)$}
&\rotatebox{90}{$(5,3^2,2^2)$}
&\rotatebox{90}{$(5^2,3,2)$}
\\\hline
(2^2,1^{11})&1&\cdot&\cdot&\cdot&\cdot&\cdot&\cdot&\cdot&\cdot&\cdot\\
(2^5,1^5)&\cdot&1&\cdot&\cdot&\cdot&\cdot&\cdot&\cdot&\cdot&\cdot\\
(3^3,2^2,1^2)&\cdot&\cdot&1&\cdot&\cdot&\cdot&\cdot&\cdot&\cdot&\cdot\\
(4,3,1^8)&q^2&q^2&\cdot&1&\cdot&\cdot&\cdot&\cdot&\cdot&\cdot\\
(4,3,2^3,1^2)&\cdot&q^2&q^2&\cdot&1&\cdot&\cdot&\cdot&\cdot&\cdot\\
(4,3^2,2,1^3)&\cdot&q^4&q^4&q^2&q^2&1&\cdot&\cdot&\cdot&\cdot\\
(4,3^2,2^2,1)&\cdot&\cdot&q^6&\cdot&q^4&q^2&1&\cdot&\cdot&\cdot\\
(5,2,1^8)&\cdot&\cdot&\cdot&q^2&\cdot&\cdot&\cdot&\cdot&\cdot&\cdot\\
(5,2^4,1^2)&\cdot&\cdot&\cdot&\cdot&q^2&\cdot&\cdot&\cdot&\cdot&\cdot\\
(5,3^2,1^4)&q^2&q^2&\cdot&q^4&q^4&q^2&\cdot&1&\cdot&\cdot\\
(5,3^2,2^2)&\cdot&q^4&\cdot&\cdot&q^6&q^4&q^2&q^2&1&\cdot\\
(5^2,1^5)&q^4&q^4&\cdot&\cdot&\cdot&\cdot&\cdot&q^2&\cdot&\cdot\\
(5^2,3,2)&q^2&q^6&\cdot&q^4&\cdot&q^2&\cdot&q^4&q^2&1\\
(5^2,4,1)&q^4&\cdot&\cdot&q^6&\cdot&q^4&\cdot&\cdot&\cdot&q^2\\
(7,3,1^5)&\cdot&\cdot&\cdot&\cdot&\cdot&q^2&\cdot&q^4&\cdot&\cdot\\
(7,3^2,2)&q^4&\cdot&\cdot&\cdot&\cdot&q^4&q^2&q^6&q^4&q^2\\
(7,6,1^2)&q^6&\cdot&\cdot&\cdot&\cdot&\cdot&\cdot&\cdot&\cdot&q^4\\
(8,2,1^5)&\cdot&\cdot&\cdot&\cdot&\cdot&q^4&\cdot&\cdot&\cdot&\cdot\\
(8,3^2,1)&\cdot&\cdot&\cdot&\cdot&\cdot&q^6&q^4&\cdot&q^2&q^4\\
(8,5,1^2)&\cdot&\cdot&\cdot&\cdot&\cdot&\cdot&\cdot&\cdot&q^4&q^6\\
(10,3,1^2)&\cdot&\cdot&\cdot&\cdot&\cdot&\cdot&q^4&\cdot&q^6&\cdot\\
(11,2,1^2)&\cdot&\cdot&\cdot&\cdot&\cdot&\cdot&q^6&\cdot&\cdot&\cdot\\
\hline
\end{array}
\]
\caption{Canonical basis in type $A^{(1)}_2$}\label{examp1}
\end{figure}

Now take $l=6$, and consider the canonical basis for the corresponding weight space in $\spin\fo$ with $h=5$. This is given by the matrix in \cref{examp2}. This weight space is spanned by $5$-strict partitions with $5$-bar-core $(12,11,7,6,2,1)$ and $5$-bar-weight $3$.

\begin{figure}[t]\footnotesize
\[
\begin{array}{c|cccccccccc|}
&\rotatebox{90}{$(12,11,7,6,5^3,2,1)$}
&\rotatebox{90}{$(12,11,10,7,6,5,2,1)$}
&\rotatebox{90}{$(15,12,11,7,6,2,1)$}
&\rotatebox{90}{$(16,12,7,6,5^2,2,1)$}
&\rotatebox{90}{$(16,12,10,7,6,2,1)$}
&\rotatebox{90}{$(16,12,11,7,5,2,1)$}
&\rotatebox{90}{$(16,12,11,7,6,2)$}
&\rotatebox{90}{$(17,12,11,6,5,2,1)$}
&\rotatebox{90}{$(17,12,11,7,6,1)$}
&\rotatebox{90}{$(17,16,11,7,2,1)$}
\\\hline
(12,11,7,6,5^3,2,1)&1&\cdot&\cdot&\cdot&\cdot&\cdot&\cdot&\cdot&\cdot&\cdot\\
(12,11,10,7,6,5,2,1)&-q^4+q^2&1&\cdot&\cdot&\cdot&\cdot&\cdot&\cdot&\cdot&\cdot\\
(15,12,11,7,6,2,1)&q^4&q^2&1&\cdot&\cdot&\cdot&\cdot&\cdot&\cdot&\cdot\\
(16,12,7,6,5^2,2,1)&q^2&q^2&\cdot&1&\cdot&\cdot&\cdot&\cdot&\cdot&\cdot\\
(16,12,10,7,6,2,1)&q^4&q^4+q^2&q^2&q^2&1&\cdot&\cdot&\cdot&\cdot&\cdot\\
(16,12,11,7,5,2,1)&\cdot&q^4&q^4&q^2&q^2&1&\cdot&\cdot&\cdot&\cdot\\
(16,12,11,7,6,2)&\cdot&\cdot&q^6&\cdot&q^4&q^2&1&\cdot&\cdot&\cdot\\
(17,11,7,6,5^2,2,1)&\cdot&\cdot&\cdot&q^2&\cdot&\cdot&\cdot&\cdot&\cdot&\cdot\\
(17,11,10,7,6,2,1)&\cdot&\cdot&\cdot&q^4&q^2&\cdot&\cdot&\cdot&\cdot&\cdot\\
(17,12,11,6,5,2,1)&q^2&q^2&\cdot&q^4&q^4&q^2&\cdot&1&\cdot&\cdot\\
(17,12,11,7,6,1)&\cdot&q^4&\cdot&\cdot&q^6&q^4&q^2&q^2&1&\cdot\\
(17,16,7,6,5,2,1)&q^4&q^4&\cdot&\cdot&\cdot&\cdot&\cdot&q^2&\cdot&\cdot\\
(17,16,11,7,2,1)&q^2&q^6&\cdot&q^4&\cdot&q^2&\cdot&q^4&q^2&1\\
(17,16,12,6,2,1)&q^4&\cdot&\cdot&q^6&\cdot&q^4&\cdot&\cdot&\cdot&q^2\\
(21,12,7,6,5,2,1)&\cdot&\cdot&\cdot&\cdot&\cdot&q^2&\cdot&q^4&\cdot&\cdot\\
(21,12,11,7,2,1)&q^4&\cdot&\cdot&\cdot&\cdot&q^4&q^2&q^6&q^4&q^2\\
(21,17,7,6,2,1)&q^6&\cdot&\cdot&\cdot&\cdot&\cdot&\cdot&\cdot&\cdot&q^4\\
(22,11,7,6,5,2,1)&\cdot&\cdot&\cdot&\cdot&\cdot&q^4&\cdot&\cdot&\cdot&\cdot\\
(22,12,11,6,2,1)&\cdot&\cdot&\cdot&\cdot&\cdot&q^6&q^4&\cdot&q^2&q^4\\
(22,16,7,6,2,1)&\cdot&\cdot&\cdot&\cdot&\cdot&\cdot&\cdot&\cdot&q^4&q^6\\
(26,12,7,6,2,1)&\cdot&\cdot&\cdot&\cdot&\cdot&\cdot&q^4&\cdot&q^6&\cdot\\
(27,11,7,6,2,1)&\cdot&\cdot&\cdot&\cdot&\cdot&\cdot&q^6&\cdot&\cdot&\cdot\\
\hline
\end{array}
\]
\caption{Canonical basis in type $A^{(2)}_4$}\label{examp2}

\end{figure}

The matrices in \cref{examp1,examp2} are almost identical, but differ because of the effect of the polynomials $K^{-1}_{\rho\pi}(-q^2)$. To see the effect of \cref{sscbv}, consider for example the fourth and fifth rows of the two matrices. In the first matrix, these rows are labelled by the partitions $\jo{(4,3,1^2)}{(1^2)}$ and $\jo{(4,3,1^2)}{(2)}$. The inverse Kostka polynomials $K^{-1}_{\rho\pi}(-q^2)$ for partitions of $2$ are given by the following matrix.
\[
\begin{array}{c|cc|}
&\rt{1^2}&\rt{2}\\\hline
(1^2)&1&\cdot\\
(2)&q^2&1\\\hline
\end{array}
\]
Hence the linear map $\sas$ sends
\begin{align*}
\jo{(4,3,1^2)}{(1^2)}\longmapsto\,&\jh{\psi(4,3,1^2)}{(1^2)}+q^2\jh{\psi(4,3,1^2)}{(2)}
\\
&=(16,12,7,6,5^2,2,1)+q^2(16,12,10,7,6,2,1)\\
\jo{(4,3,1^2)}{(2)}\longmapsto\,&\jh{\psi(4,3,1^2)}{(2)}
\\
&=(16,12,10,7,6,2,1).
\end{align*}
As a consequence, the fifth row of the matrix in \cref{examp2} is obtained from the fifth row of the matrix in \cref{examp1} by adding $q^2$ times the fourth row. A similar calculation applies to the first three rows of the matrices (in fact, the $3\times3$ matrix at the top right of the matrix in \cref{examp2} is precisely the matrix of polynomials $K^{-1}_{\rho\pi}(-q^2)$ for partitions of $3$).
\end{eg}

\section{Equivalences on weight spaces and Rouquier cores}

In the representation theory of the symmetric group, equivalences between blocks of the same defect (starting from the work of Scopes) have proved to be an important tool in studying the representation theory of symmetric groups. One of the Scopes equivalence classes of blocks was identified by Rouquier as being of particular interest. These blocks are now commonly known as \emph{Rouquier blocks} or \emph{RoCK blocks}, and (in combination with Scopes's equivalences) have been instrumental in proving various conjectures, most significantly in the proof of Brou\'e's abelian defect group conjecture for the symmetric groups \cite{ck,cr}. There is a parallel theory for blocks of double covers of symmetric groups, but with substantial additional complications. We will describe this in more detail in \cref{rocksec}, but in this section we describe analogous results in $\calf$ and $\spin\calf$: we define the notion of \emph{extremal} weight spaces in $\calf$ and $\spin\calf$, we recall the formula (due to Chuang--Tan and Leclerc--Miyachi) for the canonical bases in extremal weight spaces in type $A^{(1)}_n$, and we use this (and our main theorem from \cref{comparingsec}) to give a corresponding formula for $\spin\calf$.

\subsection{Extremal weight spaces in type $A^{(1)}_{m-1}$}

For each $m$-core $\nu$ and each $w\gs0$, let $\pw\nu w$ denote the set of partitions with $m$-core $\nu$ and $m$-weight $w$, and let $\ws\nu w$ denote the weight space in $\calf$ spanned by $\pw\nu w$. If $k\in\bbn$, we say that two weight spaces $\ws\nu w$ and $\ws\xi w$ form a \emph{$[w:k]$-pair} of residue $i$ if $\xi$ is obtained from $\nu$ by adding $k$ nodes of residue $i$. This is equivalent to saying that $\f i^k$ maps $\ws\nu w$ to $\ws\xi w$. We will write $V\cls W$ to mean that $V$ and $W$ form a $[w:k]$-pair for some $k$, and extend $\cls$ transitively to give a partial order on the set of weight spaces of $m$-weight $w$. If $V$ and $W$ form a $[w:k]$-pair of residue $i$ with $k\gs w$, we say that $V$ and $W$ are \emph{Scopes equivalent}; this is equivalent to the condition that $\f i^{(k)}$ and $\e i^{(k)}$ map the canonical basis for $V$ to the canonical basis for $W$ and vice versa. Extending this relation transitively defines an equivalence relation (the Scopes equivalence) on the set of weight spaces of $m$-weight $w$, with only finitely many equivalence classes. The partial order $\cls$ on weight spaces descends to a partial order on Scopes classes, and there is a unique maximal class for this order. We call the weight spaces in this class the \emph{extremal} weight spaces in $\calf$; these are precisely the weight spaces corresponding to RoCK blocks of symmetric groups. To describe the weight spaces in this class, we introduce the notion of a Rouquier core.

Take $w\gs1$, and let $\nu$ be an $m$-core. Say that $\nu$ is \emph{$w$-Rouquier} if it has an abacus display (with $s$ beads, say) in which there are at least $w-1$ more beads on runner $i$ than on runner $i-1$, for $i=1,\dots,m-1$. If $\nu$ is $w$-Rouquier, then $W_{\nu,w}$ is an example of an extremal weight space.

To describe the canonical basis coefficients in extremal weight spaces, we need to recall the notion of the $m$-quotient of a partition. Suppose $\nu$ is $w$-Rouquier with an $s$-bead \abd as above. By increasing $s$ by a multiple of $m$ if necessary, we can assume there are at least $w$ beads on each runner. Now given $\la\in\pw\nu w$, we can construct the \abd for $\la$ with $s$ beads from the \abd for $\nu$ by moving some beads down their runners. We define the \emph{$m$-quotient} $(\la(0),\dots,\la(m-1))$ to be the $m$-tuple of partitions in which $\la(i)$ is the partition obtained by viewing runner $i$ in isolation as $1$-runner abacus; more precisely, $\la(i)_r$ is the number of empty spaces above the $r$th lowest bead on runner $i$.

(In fact the definition of $m$-quotient we have given is not the most usual one, in which the indices $i$ may be permuted cyclically according to the number of beads on the abacus. The $m$-quotient defined here is called the \emph{ordered $m$-quotient} in \cite{mfaltredproof}.)

Now we can give the main result of Chuang--Tan and Leclerc--Miyachi (translated to the conventions we are using). Recall that $\dn\la\mu$ denotes the coefficient of $\la$ in $\cb m(\mu)$. In the following theorem, $\lr\la\si\tau$ denotes a Littlewood--Richardson coefficient, which should be interpreted as $0$ if $\card\la\neq\card\si+\card\tau$.

\begin{thm}[\bf{\cite[Theorem 1.1]{ctq}, \cite[Corollary 10]{lm}}]\label{ctmain}
Suppose $\nu$ is a $w$-Rouquier $m$-core, and $\la,\mu\in\pw\nu w$. Then $\mu$ is $m$-restricted \iff $\mu(m-1)=\vn$. If $\mu$ is $m$-restricted, then
\[
\dn\la\mu=\left(\sum\prod_{i=1}^{m-1}\lr{\la(i)}{\si(i)}{\tau(i)}\lr{\mu(i-1)}{\si(i-1)}{\tau(i)'}\right)q^{2\sum_ii(\card{\la(i)}-\card{\mu(i)})},
\]
where the sum is over all partitions $\si(1),\dots,\si(m-2),\tau(1),\dots,\tau(m-1)$, and $\si(0)$ and $\si(m-1)$ should be read as~$\la(0)$ and $\vn$ respectively.
\end{thm}

\subsection{Extremal weight spaces in type $A^{(2)}_{h-1}$}

Now we describe the corresponding theory in type $A^{(2)}_{h-1}$. If $\ga$ is an $h$-bar-core and $w\gs0$, we let $\spw\ga w$ denote the set of \bps with $h$-bar-core $\ga$ and $h$-bar-weight $w$, and we let $\sws\ga w$ denote the weight space in $\spin\calf$ spanned by $\spw\ga w$. The notion of a $[w:k]$-pair of weight spaces of bar-residue $i$, and the partial order $\cls$ are defined exactly as in type $A^{(1)}$. Two weight spaces $\sws\ga w$ and $\sws\de w$ forming a $[w:k]$-pair of bar-residue $i$ are \emph{Scopes--Kessar equivalent} if the divided powers $\f i^{(k)}$ and $\e i^{(k)}$ map the canonical basis of $\sws\ga w$ to the canonical basis of $\sws \de w$ and vice versa. (The exact condition in terms of $w,k,i$ for a Scopes--Kessar equivalence is given in \cite[Theorem 4.6]{mfspinwt2}, but we do not need this here). As in type $A^{(1)}$, we extend this relation to give an equivalence relation on the set of weight spaces with $h$-bar-weight $w$, and the order $\cls$ descends to an order on equivalence classes. Again, there is a unique maximal class, whose elements we call extremal weight spaces.

To give examples of extremal weight spaces, we introduce Rouquier bar-cores, following Kleshchev and Livesey \cite[\S4.1a]{kl}. Say that an $h$-bar-core $\ga$ is \emph{$w$-Rouquier} if in the \babd for $\ga$ there are at least $w$ beads on runner $1$, and at least $w-1$ more beads on runner $i$ than on runner $i-1$, for $2\ls i\ls n$. If $\ga$ is $w$-Rouquier, then the weight space $\spw\ga w$ is extremal.

We can use the results in this paper to find the canonical basis coefficients in extremal weight spaces. Suppose $\ga$ is a $w$-Rouquier $h$-bar-core, and let $l$ be the length of $\ga$. Then automatically $\ga\in\pnice l$ (because $\ga$ has parts congruent to $1,2,\dots,n$ modulo $h$ only), and if we apply the function $\psi$, we obtain a $w$-Rouquier $m$-core.

To give a formula analogous to \cref{ctmain}, we need to define $h$-bar-quotients. Given $\al\in\spw\ga w$, let $\al(0)$ be the partition obtained by taking all the parts of $\al$ divisible by $h$ and dividing them all by $h$. For $1\ls i\ls n$, define $\al(i)$ by letting $\al(i)_r$ be the number of empty spaces above the $r$th lowest bead on runner $i$ in the \babd for $\al$. (We remark that the particular version of $h$-bar-quotient we have defined here is specific to partitions with Rouquier cores; other versions \cite{my,yat} have been defined which apply to all $h$-strict partitions.) In particular, we can write $\al=\jh\be{\al(0)}$, where $\be\in\pnice l$, and $\jo{\psi\be}{\al(0)}$ is the partition with $m$-core $\psi\ga$ and $m$-quotient $(\al(0),\dots,\al(n))$. Moreover, the $w$-Rouquier condition implies that $\al$ is \bsep: the last bead on runner $0$ of the \babd is in position $h\al(0)_1$, while the first empty position not on runner $0$ is position $1+(d-{\al(1)}'_1)h$ or later, where $1+dh$ is the position of the last bead on runner $1$ for $\al$. Since $\al(0)_1+{\al(1)}'_1\ls w\ls d$, the \bsep condition holds.

\begin{eg}
Suppose $h=5$ and $w=4$, and let $\ga=(32,27,22,17,16,12,11,7,6,2,1)$. Then $\ga$ is a $4$-Rouquier $5$-bar-core, and $\psi\ga=(9,7,5,3^2,2^2,1^2)$ is a $4$-Rouquier $3$-core.
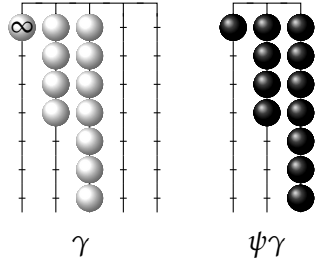
\[
\begin{array}{c@{\qquad}c}
\abacus(lmmmr,o-e\infi oonn,noonn,noonn,noonn,nnonn,nnonn,nnonn)
&
\abacus(lmr,bbb,nbb,nbb,nbb,nnb,nnb,nnb)
\\[36pt]
\ga&\psi\ga
\end{array}
\]
Now let $\al=(37,32,22,17,16,12,11,10,7,6,2,1)=\jh{(37,32,22,17,16,12,11,7,6,2,1)}{(2)}$. Then $\al$ lies in $\spw\ga4$, and has $h$-bar-quotient $((2),\vn,(1^2))$. The corresponding partition with $3$-core $\psi\ga$ and $3$-quotient $((2),\vn,(1^2))$ is $(12,10,5,3^2,2^5,1^2)$.
\[
\begin{array}{c@{\qquad}c}
\abacus(lmmmr,o-e\infi oonn,noonn,ooonn,noonn,nnonn,nnnnn,nnonn,nnonn)
&\abacus(lmr,nbb,nbb,bbb,nbb,nnb,nnn,nnb,nnb)
\\[43pt]
\al&(12,10,5,3^2,2^5,1^2)
\end{array}
\]
\end{eg}
Now we can apply \cref{sscbv,ctmain}, and we obtain the following result giving the canonical bases for weight spaces with Rouquier bar-cores. (This generalises in a straightforward way to give the canonical bases for all extremal weight spaces.)

\begin{thm}\label{mainrouq}
Suppose $\ga$ is a $w$-Rouquier $h$-bar-core, and that $\al$ and $\be$ are \bps with $h$-bar-core $\ga$ and $h$-bar-weight $w$. Then $\be$ is restricted \iff $\be(n)=\vn$. If $\be$ is restricted, then
\[
\sdn\al\be=\left(\sum K^{-1}_{\al(0)\si(0)}(-q^2)\prod_{i=1}^n\lr{\al(i)}{\si(i)}{\tau(i)}\lr{\be(i-1)}{\si(i-1)}{\tau(i)'}\right)q^{2\sum_ii(\card{\al(i)}-\card{\be(i)})},
\]
where the sum is over all partitions $\si(0),\dots,\si(n-1),\tau(1),\dots,\tau(n)$, and $\si(n)$ should be read as~$\vn$.
\end{thm}

\begin{eg}
Take $h=3$ and $w=4$. The $3$-bar-core $(10,7,4,1)$ is $4$-Rouquier, and the canonical basis in the corresponding weight space of $3$-bar-weight $4$ is given by the matrix in \cref{examp3}.
\begin{figure}[ht]
\[
\begin{array}{c|ccccc|}
&\rotatebox{90}{$(10,7,4,3^4,1)$}
&\rotatebox{90}{$(10,7,6,4,3^2,1)$}
&\rotatebox{90}{$(10,7,6^2,4,1)$}
&\rotatebox{90}{$(10,9,7,4,3,1)$}
&\rotatebox{90}{$(12,10,7,4,1)$}
\\\hline
(10,7,4,3^4,1)&1&\cdot&\cdot&\cdot&\cdot\\
(10,7,6,4,3^2,1)&q^6-q^4+q^2&1&\cdot&\cdot&\cdot\\
(10,7,6^2,4,1)&-q^6&q^2&1&\cdot&\cdot\\
(10,9,7,4,3,1)&-q^6+q^4&q^2&q^2&1&\cdot\\
(12,10,7,4,1)&q^6&q^4&\cdot&q^2&1\\
(13,7,4,3^3,1)&q^2&q^2&\cdot&\cdot&\cdot\\
(13,7,6,4,3,1)&-q^6+q^4&-q^6+q^4+q^2&q^2&q^2&\cdot\\
(13,9,7,4,1)&q^6&q^6+q^4&q^4&q^4+q^2&q^2\\
(13,10,4,3^2,1)&\cdot&q^4&\cdot&q^4&\cdot\\
(13,10,6,4,1)&\cdot&q^6&q^4&q^6+q^4&q^4\\
(16,7,4,3^2,1)&q^4&q^4&q^4&\cdot&\cdot\\
(16,7,6,4,1)&q^6&q^6+q^4&q^6&q^4&\cdot\\
(13,10,7,3,1)&\cdot&\cdot&\cdot&q^6&q^6\\
(16,10,4,3,1)&\cdot&q^6&q^6&q^6&\cdot\\
(19,7,4,3,1)&q^6&q^6&\cdot&\cdot&\cdot\\
(13,10,7,4)&\cdot&\cdot&\cdot&\cdot&q^8\\
(16,10,7,1)&\cdot&\cdot&\cdot&q^8&\cdot\\
(16,13,4,1)&\cdot&\cdot&q^8&\cdot&\cdot\\
(19,10,4,1)&\cdot&q^8&\cdot&\cdot&\cdot\\
(22,7,4,1)&q^8&\cdot&\cdot&\cdot&\cdot\\\hline
\end{array}
\]
\caption{Canonical basis in an extremal weight space}\label{examp3}
\end{figure}
As an example of how entries are calculated, consider the $(\al,\be)$-entry where $\al=(13,7,6,4,3,1)$ and $\be=(10,7,6,4,3^2,1)$. In this case
\[
(\al(0),\al(1))=((2,1),(1)),\qquad(\be(0),\be(1))=((2,1^2),\vn),
\]
giving $\sum_ii(\card{\al(i)}-\card{\be(i)})=1$. The Littlewood--Richardson coefficient $\lr{(2,1^2)}\si{(1)}$ equals $1$ if $\si=(1^3)$ or $(2,1)$, and $0$ otherwise. So
\[
\sdn\al\be=(K^{-1}_{(2,1)(1^3)}(-q^2)+K^{-1}_{(2,1)(2,1)}(-q^2))q^2=((q^2-q^4)+1)q^2.
\]
\end{eg}

\section{RoCK blocks of symmetric groups and their double covers}\label{rocksec}

We end this paper by explaining the connections between our results and the spin representation theory of the symmetric group, giving conjectural decomposition numbers for RoCK blocks with abelian defect.

We refer to the book by Kleshchev \cite{kleshbook} for background on spin representation theory of $\sss n$, but we summarise the main points here. As in the introduction, let $\hsss n$ denote a proper double cover of $\sss n$. We consider representations of $\hsss n$ over an algebraically closed field $\bbf$ of odd characteristic $h$. The non-trivial central element $z\in\hsss n$ acts as $\pm1$ on any irreducible module; modules on which $z$ acts as $-1$ are called \emph{spin modules}. A block of $\hsss n$ is called a spin block if it contains a spin module (in which case all the irreducible modules in the block are spin modules).

In practice it is more convenient to work with $\bbf\hsss n$ as a superalgebra, and consider supermodules. So we will consider spin superblocks (which coincide with spin blocks except in the trivial case of blocks of defect $0$). For each strict partition $\al$ of $n$ (i.e.\ each partition in which the positive parts are distinct) there is an irreducible spin supermodule for $\bbc\hsss n$, and these modules give all the irreducible spin supermodules for $\bbc\hsss n$. We let $\sspe\al$ denote an $h$-modular reduction of this module. For each restricted $h$-strict partition $\be$ of $n$, there is an irreducible spin supermodule $\sid\be$ for $\bbf\hsss n$, and these modules give all the irreducible spin supermodules for $\bbf\hsss n$. So (apart from the distinction between modules and supermodules) the decomposition number problem for spin representations of $\hsss n$ asks for the composition multiplicities $[\sspe\al:\sid\be]$ for all $\al,\be$.

The block classification for $\hsss n$ (due to Humphreys \cite{hum}) says that $\sspe\al$ and $\sid\be$ lie in the same superblock of $\hsss n$ \iff $\al$ and $\be$ have the same $h$-bar-core (in which case they have the same $h$-bar-weight as well). In particular, we can talk about the $h$-bar-core and $h$-bar-weight of a superblock, and there is a direct correspondence between spin superblocks of $\hsss n$ (as $n$ varies) and weight spaces in the Fock space $\spin\calf$. Moreover, the action of the generators $\spe i,\spf i$ on standard basis elements in $\spin\calf$ corresponds to the branching rules describing induction and restriction of the modules $\sspe\al$. This led Leclerc and Thibon to draw connections between canonical basis coefficients and decomposition numbers. To give a statement of their conjecture, we introduce some more notation. Given an $h$-strict partition $\al$ of $n$, we can define in a combinatorial way a restricted $h$-strict partition called the \emph{regularisation} $\al\reg$ of $\al$; this was introduced by Brundan and Kleshchev, who showed that the decomposition number $[\sspe\al:\sid{\al\reg}]$ is non-zero (giving its value explicitly), and that $\al\reg$ is the most dominant partition with this property. Now given any restricted $h$-strict partition $\be$ of $n$, define the \emph{divided decomposition number}
\[
D_{\al\be}=\frac{[\sspe\al:\sid\be]}{[\sspe\al:\sid{\al\reg}]}.
\]

Then the Leclerc--Thibon conjecture on spin decomposition numbers states that $D_{\al\be}$ is simply the evaluation of $\sdn\al\be$ at $q=1$ when $n<h^2$. A natural extension of this conjecture (analogous to the James conjecture for the symmetric groups) would weaken the condition $n<p^2$ to include all cases where $\al$ and $\be$ have $h$-bar-weight less than $h$ (this condition corresponds to the block containing $\sspe\al$ having abelian defect groups). As mentioned in the introduction, this conjecture is known to be false (in fact it predicts negative decomposition numbers!) but it is nevertheless true in many cases; in particular, in all blocks of abelian defect for which the decomposition numbers are explicitly known, the Leclerc--Thibon conjecture holds.

Now define a \emph{RoCK block} to be a superblock with $h$-bar-weight $w$ whose $h$-bar-core is $w$-Rouquier. Applying \cref{mainrouq}, we make the following explicit conjecture. Here we write $K^{-1}_{\al\be}$ for the specialisation of $K^{-1}_{\al\be}(t)$ at $t=-1$. We continue to write $n=\frac12(h-1)$.

\begin{conj}\label{rockconj}
Suppose $\ga$ is a $w$-Rouquier $h$-bar-core, and that $\al$ is a strict partition and $\be$ a restricted $h$-strict partition, both with $h$-bar-core $\ga$ and $h$-bar-weight $w<h$. Define the $h$-bar-quotients $(\al(0),\dots,\al(n))$ and $(\be(0),\dots,\be(n))$ as above. Then
\[
D_{\al\be}=\sum K^{-1}_{\al(0)\si(0)}\prod_{i=1}^n\lr{\al(i)}{\si(i)}{\tau(i)}\lr{\be(i-1)}{\si(i-1)}{\tau(i)'}.
\]
where the sum is over all partitions $\si(0),\dots,\si(n-1),\tau(1),\dots,\tau(n)$, and $\si(n)$ should be read as~$\vn$.
\end{conj}

Recent work of Kleshchev and Livesey \cite{kl} studies RoCK blocks in detail, particularly in the abelian defect case. In a forthcoming paper \cite{fkm}, we will use these results and the results of the present paper to prove \cref{rockconj}.

\end{document}